\documentclass[12pt]{amsart}
\newfont{\sheaf}{eusm10 scaled\magstep1}

\usepackage{amssymb,amsbsy,amscd}

\input xy
\xyoption{all}
\CompileMatrices

\newcommand{\ra}{\ensuremath{\rightarrow}}

\def\eea{\end{eqnarray*}}
\def\bea{\begin{eqnarray*}}

\newcommand{\Proof}{{\it Proof. }}
\newcommand{\QED}{{\hfill $Q.E.D.$}}

\newtheorem{teo}{Theorem}[section]
\newtheorem{df}[teo]{Definition}
\newtheorem{lem}[teo]{Lemma}
\newtheorem{cor}[teo]{Corollary}

\newtheorem{oss}[teo]{Remark}
\newtheorem{prop}[teo]{Proposition}

\newcommand{\C}{\ensuremath{\mathbb{C}}}

\newcommand{\Z}{\ensuremath{\mathbb{Z}}}

\newcommand{\F}{\ensuremath{\mathbb{F}}}
\newcommand{\M}{\ensuremath{\mathbb{M}}}
\newcommand{\N}{\ensuremath{\mathbb{N}}}

\newcommand{\BBB}{\ensuremath{\mathbb{U}}}
\newcommand{\TT}{\ensuremath{\mathbb{T}}}
\newcommand{\PP}{\ensuremath{\mathbb{P}}}

\newcommand{\SSS}{\ensuremath{\mathcal{S}}}
\newcommand{\AAA}{\ensuremath{\mathcal{A}}}

\begin{document}

\title{ Beauville surfaces without real structures, I.}

\author{Ingrid Bauer- Fabrizio Catanese - Fritz Grunewald}\
   \footnote{  Universit\"at Bayreuth$^2$- Universit\"at D\"usseldorf\\
The
research of the  authors was performed in the realm  of the
   SCHWERPUNKT "Globale Methode in der komplexen Geometrie", and of the
EAGER EEC Project.  }
   \\

\date{August 2 , 2004}
\maketitle

\begin{abstract}
Inspired by a construction by Arnaud Beauville of a surface
of general type with $K^2 = 8, \ p_g =0$,  the second author defined
the Beauville
surfaces  as the surfaces which are rigid, i.e., they have no nontrivial
deformation, and admit un unramified  covering
which is  isomorphic to a product of curves of genus at least $2$.

In this case the moduli space of surfaces homeomorphic to the given surface
consists either of a unique real point, or of a pair of complex conjugate
points corresponding to complex conjugate surfaces. It may also 
happen that a Beauville
surface is biholomorphic to its complex conjugate surface, neverless 
it fails to admit a
real structure.

First aim of this note is to provide series of concrete examples
   of the second situation, respectively of the third.

Second aim is to introduce a wider audience, especially group theorists,
   to the problem of classification of such surfaces, especially  with regard
to the  problem of existence of real structures on them.
    \end{abstract}


\tableofcontents

\section{Introduction}

In \cite{bea} (see p. 159) A. Beauville constructed a new surface
of general type with $K^2 = 8, \ p_g =0$ as a quotient of the
product of two Fermat curves of degree $5$ by the action of the group
$(\Z/5\Z)^2$. Inspired by this construction, in the
article \cite{cat00}, dedicated
to the geometrical properties of varieties which admit an unramified covering
biholomorphic to a product of curves, the following definition was given

\begin{df}
A Beauville surface is a compact complex surface $S$ which

1) is rigid , i.e., it has no nontrivial deformation,

2) is isogenous to a higher product, i.e., it
   admits an unramified covering which is isomorphic (i.e., biholomorphic) to
the product of two curves $C_1$, $C_2$ of genera $\geq 2$.
\end{df}

It was proven in \cite{cat00} (cf. also \cite{cat03}) that any
surface $S$ isogenous to a higher product has a unique minimal realization
   as a quotient $S = (C_1 \times C_2)/G$,
where $G$ is a finite group acting freely and with the property that no
element acts trivially on one of the two factors $C_i$. Moreover,
any other smooth surface $X$ with the same topological Euler number
as $S$, and with isomorphic fundamental group, is diffeomorphic to $S$,
and either $X$ or the conjugate surface $\bar{X}$ belongs to
an irreducible family of surfaces containing $S$ as element.

Therefore, in the case that $S$ is a Beauville surface, either
$X$ is isomorphic to $S$, or $X$ is isomorphic to $\bar{S}$.

In order to reduce the description of Beauville surfaces to some
group theoretic statement, we need to recall that surfaces
isogenous to a higher product belong to two types:

\begin{itemize}
\item
$S$ is of {\bf unmixed type} if the action of $G$ does not mix the two factors,
i.e., it is the product action of respective actions of $G$ on $C_1$,
resp. $C_2$. We set then $G^0:=G$.
\item
$S$ is of {\bf mixed type}, i.e., $C_1$ is isomorphic to $C_2$, and the
subgroup $G^0$ of transformations in $G$ which do not mix the factors
has index precisely $2$ in $G$.
\end{itemize}

It is obvious from the above definition that every Beauville surface 
of mixed type has
an unramified double covering which is a Beauville surface of unmixed type.

The rigidity property of the Beauville
surface is equivalent to the fact that $C_i /G \cong \PP^1$ and that
the projection $C_i \ra C_i /G \cong \PP^1$ is branched in three points.
Therefore the datum of a Beauville surface of unmixed type
is determined, once we look at the monodromy of each covering of
$\PP^1$, by the datum
of a finite group
$G = G^0 $ together with two respective systems of generators, $(a,c)$ and
$(a',c')$, which satisfy a further property (*), ensuring that
the product action of $G$ on $ C \times C'$ is free, where
$C := C_1$, $C' : = C_2$ are the corresponding curves  with an action of $G$
associated to the monodromies determined by $(a,c)$, resp.
$(a',c')$.

Define $b,\, b'$ by the properties $abc= a'b'c' = 1$,
let $\Sigma$  be the union of the conjugates of the cyclic
subgroups generated by $a,b,c$ respectively, and
define $\Sigma'$
analogously: then property (*) is the following
$$ (*) \  \Sigma \cap \Sigma' = \{ 1_G \} .  $$
In the mixed case, one requires instead that the two systems
of generators be related by an automorphism $\phi$ of $G^0$
which should satisfy the further conditions:

\begin{itemize}
\item
$\phi^2$ is an inner automorphism, i.e., there is an element
   $\tau \in G^0$ such that $\phi^2 = {\rm Int}_{\tau}$,
\item
$\Sigma \cap \phi (\Sigma) = \{ 1_{G^0} \}, $
\item
there is no $g \in G^0$ such that $ \phi(g) \tau g \in \Sigma$,
\item
moreover $\phi (\tau) = \tau$ and indeed the elements in the trivial coset
of $G^0$  are transformations of $C \times C$ of the form
$$ g(x,y) = ( g (x) , \phi (g) (y) $$
while transformations in the nontrivial coset are transformations of the
form
$$ \tau' g (x,y) = ( \phi (g) (y), \tau g (x) ) .$$
\end{itemize}

\begin{oss}
The choice of $\tau$ is not unique, we can always replace $\tau$ by
   $ \phi(g) \tau g $, where $g \in G^0$ is arbitrary,  and accordingly replace
   $ \phi $ by $ \phi \circ  {\rm Int}_g$.
\end{oss}

\bigskip

Observe however that, in the case where $G^0$ has only inner automorphisms,
certainly we cannot find any Beauville surface of mixed type,
since the second of the above properties will obviously be violated.
\bigskip

In this paper we use the above definition of Beauville surfaces of
unmixed and mixed
type to create group-theoretic data which will allow us to treat the
following problems:

\bigskip
\centerline{\it 1. The biholomorphism
problem for Beauville surfaces}

We introduce sets of structures $\BBB(G)$ and $\M(G)$ for every finite
group $G$
together with groups ${\rm A}_\BBB(G)$, ${\rm A}_\M(G)$  acting on  them.
We call $\BBB(G)$ the set of unmixed Beauville
structures and $\M(G)$ the set of
mixed Beauville structures on $G$. Using constructions from
\cite{cat00} and \cite{cat03} we associate an unmixed Beauville surface
$S(v)$ to every $v\in\BBB(G)$ and a mixed Beauville surface $S(u)$
to every $u\in\M(G)$. The minimal Galois representation of
every Beauville surface yields a surface $S(v)$ in
the unmixed case, respectively a surface $S(u)$ in the mixed case. 
We then show
that
$S(v)$ is biholomorphic to $S(v')$ ($v,v'\in\BBB(G)$) if and only if
$v$ lies in the  ${\rm A}_\BBB(G)$ orbit of $v'$. We also establish the
analogous result in the mixed case.

\bigskip
\centerline{\it 2. Existence and classification problem for Beauville surfaces}

The existence problem now asks for finite groups $G$ such that $\BBB(G)$ or
$\M(G)$ is not empty. So far only abelian groups $G$ were known with
$\BBB(G)\ne \emptyset$. We give many more examples of groups $G$ with
$\BBB(G)\ne \emptyset$.
In the mixed case it is not
immediately clear that the requirements for the corresponding structures
can be met. In fact no examples were known previously.
We give here a group theoretic construction which produces
finite groups $G$ with $\M(G)\ne\emptyset$.

The classification problem has two meanings. First of all we might like to
find all groups $G$ with $\BBB(G)\ne \emptyset$ or $\M(G)\ne \emptyset$.
In \cite{cat03} all finite abelian groups $G$ are found with
$\BBB(G)\ne \emptyset$ (we give a proof of this fact in Section 3, Theorem
\ref{abelian}). We show here amongst other things that a group
$G$ with $\BBB(G)\ne \emptyset$ cannot be a non-trivial quotient of one of the
non hyperbolic triangle groups.

Our examples however show that this classification problem might be hopeless.
In fact we show in Section 3.2 that every 
finite group $G$ of exponent $n$ with $G.C.D.(n,6)=1$, which
is generated by $2$ elements and which has $\Z/n\Z\times \Z/n\Z$ as
abelianization  has $\BBB(G)\ne \emptyset$. Even if, by Zelmanov's
solution of the restricted Burnside problem, there is for every $n$ a
maximal such group, the number of groups involved here is vast. We also show
  
\begin{teo}\label{intro1}
Let $G$ be one of the groups
${\bf SL}(2,\F_p)$ or ${\bf PSL}(2,\F_p)$ where $\F_p$ is the prime field with
$p$ elements and where the prime $p$ is distinct from $2,3,5$. Then there is
an unmixed Beauville surface with group $G$.
\end{teo}

A finer classification entails the determination of all orbits of Beauville
structures for a fixed group $G$ or for an interesting series of groups.
  We do not address this problem here.

\medskip
\centerline{\it 3. Is $S$  biholomorphic to $\bar{S}$? }

We give here examples of Beauville surfaces $S$
such that the complex conjugate surface $\bar{S}$ is not biholomorphic to $S$.
Note that $\bar{S}$ is
the same differentiable manifold as $S$,  but with complex structure
$-J$ instead of $J$. To do this we introduce involutions
$\iota: \BBB(G)\to\BBB(G)$ and $\iota: \M(G)\to\M(G)$. We prove that
$S(v)$ ($v\in\BBB(G)$) is biholomorphic to
$\overline{S(v)}$ if and only if $v$ is in the
${\rm A}_\BBB(G)$ orbit of $\iota(v)$. We also show the analogous result in
the mixed case. We use this to produce the following explicit example:

\begin{teo}\label{sym} Let $G$ be the symmetric group $\SSS_n$ in $n
\geq 8$ letters,
let $S (n)$  be the unmixed Beauville surface corresponding to the choice of
$a:=(5,4,1) (2,6)$, $ c:=(1,2,3)(4,5,.\dots, n)$,
and of $ a' : =  \sigma^{-1}$, $c' : = \tau \sigma^2$, where $\tau := (1,2)$
and $\sigma: = (1,2, \dots , n)$.

Then $S(n)$ is not biholomorhic to $\overline{S(n)}$ provided that
$n \equiv 2 \ {\rm mod} \ 3$.
\end{teo}

\medskip
We shall now give the construction of a mixed Beauville surface with the same
property. We shall first describe the group $G$ and its subgroup $G^0$.  

Let $H$ be nontrivial group. Let
$\Theta :H\times H \to H\times H$ be the automorphism defined by
$\Theta(g,h):=(h,g)$ ($g,h\in H$). We consider the semidirect product
\begin{equation}
G:=H_{[4]}:= (H\times H) \rtimes  \Z / 4\Z
\end{equation}
where the generator $1$ of $ \Z / 4\Z$ acts through $\Theta$ on $H\times H$.
Since $\Theta^2$ is the identity we find
\begin{equation}
G^0:=H_{[2]}:=H\times H\times 2\Z / 4\Z\cong H\times H\times \Z / 2\Z
\end{equation}
as a subgroup of index $2$ in $H_{[4]}$.
\begin{teo}\label{intro2}
Let $p$ be a prime with $p\equiv 3$ mod $4$ and 
 $p\equiv 1$ mod $5$ and consider the group $H:={\bf SL}(2,\F_p)$. Let $S$ be
 the mixed Beauville surface corresponding to the data $G:=H_{[4]}$, 
$G^0:=H_{[2]}$ and to a certain system of generators $(a,c)$ of $H_{[2]}$ with 
${\rm ord}(a)=20$, ${\rm ord}(c)=30$, ${\rm ord}(a^{-1}c^{-1})=5p$.

Then $S$ is not biholomorhic to $\overline{S}.$
\end{teo}

Different types of examples of rigid surfaces $S$ which are not
isomorphic to $\bar{S}$ have been constructed in
   \cite{k-k02}, using Hirzebruch type examples of ball quotients.

\medskip
\centerline{\it 4. Is $S$  real? }

A surface $S$ is said to be real if there exists a biholomorphism 
$\sigma:S\to \bar S$
with the
property that $\sigma^2 = {\rm Id}$. In this case we say that $S$ has a real
structure. We translate this problem into group theory and obtain the
following examples.

\begin{teo}\label{intro3}
Let $p>5$ be a prime with $p\equiv 1$ mod $4$, $p\not\equiv 2,4$ mod 5,  
$p\not\equiv 5$ mod $13$ and $p\not\equiv 4$ mod $11$. Set $n:=3p+1$. 
Then
there is  an unmixed Beauville surface $S$ with group $\AAA_n$
which is biholomorphic to the complex conjugate surface $\bar{S}$ but is not
real.
\end{teo}

Further examples of real and non real Beauville surfaces will be 
given in the sequel to
this paper.

\medskip
{\it Acknowledgement:} We thank Benjamin Klopsch for help with the alternating
groups.

\bigskip

\section{Triangular curves and group actions}\label{triang}

In this section we recall the construction of triangular curves as given in
\cite{cat00}, \cite{cat03}. They are the building blocks for the
Beauville surfaces of both unmixed and mixed type. We add some
group-theoretic observations which will help with the  classification problems
of Beauville surfaces which were mentioned above and which will be 
studied later.

We need the following group theoretic notation. Let $G$ be a group and
let $M,\, N$ be two sets equipped with an action (from the left) of
$G$. We say that a map
$\sigma :M \to N$ is {\it $G$-twisted-equivariant} if there is an automorphism
$\psi: G\to G$ of $G$ with
\begin{equation}
\sigma(g P)=\psi(g)\sigma(P) \quad {\rm for\ all} \ g\in G, P\in M.
\end{equation}

Let now $G$ be a finite group and
$(a,c)$ a pair of elements of $G$. We define
\begin{equation} \Sigma(a,c):=\bigcup_{g\in G}
\bigcup_{i=0}^\infty\ \{ga^ig^{-1},gc^ig^{-1},g(ac)^ig^{-1}\} \end{equation}
to be the union of the $G$-conjugates of the cyclic groups generated by $a$,
$c$ and $ac$. Moreover set
\begin{equation}
\mu (a,c) := \frac{1}{{\rm ord}(a)} + \frac{1}{{\rm ord}(c)}
+ \frac{1}{{\rm ord}(ac)},
\end{equation}
where ${\rm ord}(a)$ stands for the order of the element $a\in G$. We
furthermore call
\begin{equation} ({\rm ord}(a),{\rm ord}(c),{\rm ord}(ac))  \end{equation}
the {\it type} of the pair $(a,c)$ and define
\begin{equation} \nu(a,c):={\rm ord}(a){\rm ord}(c){\rm ord}(ac).
\end{equation}
We consider here finite groups $G$ having a pair $(a,c)$ of generators. Setting
$(r,s,t):= ({\rm ord}(a),{\rm ord}(c),{\rm ord}(ac))$, such a group is a
quotient of the triangle group
\begin{equation} T(r,s,t):=\langle x,y \ |\ x^r=y^s=(xy)^t=1\rangle.
\end{equation}
We define
\begin{equation}
\TT(G):=\{ (a,c)\in G\times G\ |\ \langle a,c\rangle =G\,\}.
\end{equation}
Given $(a,c)\in \TT(G)$ we shall consider the {\it triangular triple}
$(a,b,c):=(a,a^{-1}c^{-1},c)$ attached to $(a,c)$.
Clearly the automorphism group ${\rm Aut}(G)$ of $G$ acts diagonally on
$\TT(G)$. If $\TT(G)\ne \emptyset $ then this action is faithful.
We define additionally the following permutations of $\TT(G)$:
\begin{equation}\label{si1}
\sigma_0 :(a,c)\mapsto (a,c),\
\sigma_1 :(a,c)\mapsto (a^{-1}c^{-1},a),\
\sigma_2 :(a,c)\mapsto (c,a^{-1}c^{-1}),\
\end{equation}
\begin{equation}\label{si2}
\sigma_3 :(a,c)\mapsto (c,a),\
\sigma_4 :(a,c)\mapsto (c^{-1}a^{-1},c),\
\sigma_5 :(a,c)\mapsto (a,c^{-1}a^{-1}).
\end{equation}

Observe that the set  $\TT(G)$ is in bijection with the set of triples
$\TT_{\rm tr}(G): = \{ (a,b,c) | abc = 1\}$.
By looking at those triples we see that $\sigma_0$ is the identity,
$\sigma_1$ is the 3-cycle $(a,b,c) \ra (b,c,a)$,
$\sigma_3$ is the permutation $(a,b,c) \ra (c,c^{-1}bc,a)$, while 
$\sigma_2 = \sigma_1
^2$ and $\sigma_1 \sigma_3 = \sigma_4,\  \sigma_1^2 \sigma_3 = \sigma_5$.
We see therefore that we have the relations:
\begin{equation}\label{si3}
\sigma_1^3= \sigma_3^2 = \sigma_0,\  \sigma_2 = \sigma_1 ^2,\
\sigma_1 \sigma_3 = \sigma_4,\  \sigma_1^2 \sigma_3 = \sigma_5,\
\end{equation}
\begin{equation}\label{si4}
( \sigma_1 \sigma_3  ) ^2 = \sigma_4^2={\rm
Int}_{c^{-1}}\circ \sigma_0.
\end{equation}
Let us write
\begin{equation}
{\rm A}_\TT(G):=\langle {\rm Aut}(G), \sigma_1,\ldots ,\sigma_5\rangle
\end{equation}
for the permutation group generated by these operations. The above 
equations show that
we have a homomorphism of the symmetric group $\SSS_3$ into ${\rm 
A}_\TT(G) / {\rm Int
(G)}$ and that
${\rm Aut}(G)$ is a normal subgroup of index $\le 6$ in ${\rm A}_\TT(G)$, with
quotient a subgroup of $\SSS_3$ .
In particular, every
element $\rho\in {\rm A}_\TT(G)$ can be written as
\begin{equation}
\rho=\psi\circ\sigma_i
\end{equation}
for an automorphism $\psi$ of $G$ and an element $\sigma_i$ from the 
above list. 

We also define 
\begin{equation}
{\rm I}_\TT(G):=\langle {\rm Int}(G), \sigma_1,\ldots ,\sigma_5\rangle
\end{equation}
where ${\rm Int}(G)$ the (normal) subgroup of
${\rm A}_\TT(G)$ consisting of the inner automorphisms.

By an operation from ${\rm A}_\TT(G)$ we may always ensure that a pair
$(a,c)\in\TT(G)$ satisfies
$${\rm ord}(a)\le {\rm ord}(b)={\rm ord}(a^{-1}c^{-1})\le {\rm ord}(c)$$
in which case we call the pair {\it normalised}. We call $(a,c)$ {\it strict},
if all inequalities are strict, {\it critical} if all the three orders are
equal, and  {\it subcritical}  otherwise.

We shall attach now a complex curve $C(a,c)$ to every pair $(a,c)\in
\TT(G)$. It will be constructed as a ramified covering of $\PP^1_\C$.
Consider  the set $B \subset \PP^1_{\C}$ consisting of  three
real points $ B : = \{-1, 0, 1\}.$
We choose $\infty$ as a base point in $\PP^1_{\C} - B$, and we take
the following generators $ \alpha, \beta, \gamma$ of
$ \pi_1 (\PP^1_{\C} - B, \infty)$ :

\begin{itemize}
\item
$\alpha$ goes from $\infty$ to $-1 -\epsilon$
along the real line, passing through $ -2$, then makes
a full turn
counterclockwise around the circumference with centre $-1$ and radius
$\epsilon$, then goes back to $2$ along the same way on the real line.
\item
$\gamma$ goes from $\infty$ to $1 +\epsilon$
along the real line,  then makes a full turn
counterclockwise around the circumference with centre $+1$ and radius
$\epsilon$, then goes back to $\infty$ along the same way on the real line.
\item
$\beta$ goes from $\infty$ to $1 +\epsilon$
along the real line,   makes a half turn
counterclockwise around the circumference with centre $+1$ and radius
$\epsilon$, reaching $1 -\epsilon$, then proceeds
along the real line  reaching $+ \epsilon$,
   makes a full turn
counterclockwise around the circumference with centre $0$ and radius
$\epsilon$,  goes back to $1 -\epsilon$ along the same way
on the real line, makes again  a half turn
clockwise around the circumference with centre $+1$ and radius
$\epsilon$, reaching $1 +\epsilon$, finally it proceeds
along the real line returning to $\infty$.
\end{itemize}

A graphical picture of  $\alpha,\, \beta,\, $ is:

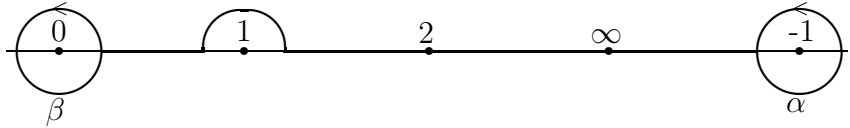
\font\mate=cmmi8
\def\min{\hbox{\mate \char60}\hskip1pt}
\begin{center}
\begin{picture}(350,100)(0,-40)
%
%
\put(0,0){\line(1,0){300}}
\put(238,0){\line(1,0){80}}
%
%
\put(20,0){\circle*{3}}
\put(17,4){0}
\put(90,0){\circle*{3}}
\put(87,4){1}
\put(160,0){\circle*{3}}
\put(156,3){2}
\put(228,0){\circle*{3}}
\put(221,3){$\infty$}
\put(300,0){\circle*{3}}
\put(296,4){-1}
{\thicklines
%
%
\put(20,0){\circle{30}}
\put(36,0){\line(1,0){39}}
\put(90,0){\oval(31,31)[t]}
\put(105,0){\line(1,0){55}}
\put(17,14){$\min$}
\put(15,-26){$\beta$}
%
%
\put(300,0){\circle{30}}
\put(160,0){\line(1,0){124}}
\put(297,14){$\min$}
\put(295,-23){$\alpha$}
}
\end{picture}
\end{center}

Writing $ \alpha,  \beta, \gamma$ for the corresponding elements of
  $ \pi_1 (\PP^1_{\C} - B, \infty)$ we find
$$\pi_1 (\PP^1_{\C} - B, \infty)=\langle \alpha, \beta, \gamma \ | \
\alpha \beta \gamma = 1\rangle$$
and $ \alpha,  \gamma$ are  free generators
of $ \pi_1 (\PP^1_{\C} - B, \infty)$.

Let now $G$ be a finite group and $(a,c)\in \TT(G)$.
Observe that by Riemann's existence theorem the elements $a$,
$b=a^{-1}c^{-1}$, $c$, once we fixed a basis of the fundamental group
of $\PP^1_\C - \{-1,0,1\}$ as above, give rise to a surjective
homomorphism
\begin{equation}\label{e1}
\pi_1 (\PP^1_{\C} - B, \infty)\to G,\qquad \alpha\mapsto a, \
\gamma \mapsto c
\end{equation}
and to a Galois covering
$\lambda:C \rightarrow \PP^1_\C$
ramified only in $\{-1,0,1\}$ with ramification indices equal to the
orders of $a$, $b$, $c$ and with group $G$ (beware, this means that these data
yield a well determined action of $G$ on $C(a,c)$).

We call this covering a  {\it
triangular covering}. We embed
$G$ into $\SSS_{G}$ as the transitive subgroup  of left translations.
The monodromy homomorphism \footnote{Actually, with the usual 
conventions the monodromy
is an antihomomorphism; there are two ways to remedy this problem, 
here we shall
do it by considering the composition of paths $ \gamma \circ \delta$ 
as the path obtained
by following first
$\delta$ and then
$\gamma$ . }
$$m_\lambda: \pi_1 (\PP^1_{\C} - B, \infty)\to \SSS_{G}$$
maps onto the embedded
subgroup $G$ and is the same as the homomorphism (\ref{e1}).

Note that by
Hurwitz's formula we have
\begin{equation}\label{gen}
g(C(a,c)) = 1 + \frac{1- \mu (a,c)}{2} \ |G|
\end{equation}
for the genus $g(C(a,c))$ of the curve $C(a,c)$.

Let now $(a,c),\, (a',c')\in\TT(G)$. A {\it twisted covering 
isomorphism} from the
Galois covering $\lambda:C(a,c) \rightarrow \PP^1_\C$ to the Galois covering
$\lambda':C(a',c') \rightarrow \PP^1_\C$ is a pair $(\sigma,\delta)$
of biholomorphic maps $\sigma : C(a,c)\to C(a',c')$ and
$\delta :\PP^1_\C\to\PP^1_\C$ with $\delta(B)=B$ such that the diagram
\begin{equation}\def\normalbaselines{\baselineskip20pt
   \lineskip3pt\lineskiplimit3pt }
\def\mapright#1{\smash{
      \mathop{\longrightarrow}\limits^{#1}}}
\def\mapdown#1{\Big\downarrow
\rlap{$\vcenter{\hbox{$\scriptstyle#1$}}$}}
\begin{matrix}
&C(a,c)&\mapright{\sigma}&
C(a',c')&\cr
&{\bigg\downarrow} {\lambda}&&{\bigg\downarrow} {\lambda'}\cr
&\PP^1_\C
&\mapright{\delta}&
\PP^1_\C.
&\cr\end{matrix}
\end{equation}
is commutative. We shall moreover say that we have a {\it strict covering 
isomorphism} if moreover the
map $\delta$ is the identity.

Consider now $G$ as acting on $C(a,c)$ as a group of covering 
transformations over
$\lambda$, and conjugate a transformation $ g \in G$ by $\sigma$: 
since  $\sigma \circ g
\circ
\sigma^{-1}$ is a covering transformation of $C(a',c')$, we obtain in 
this way an
automorphism
$\psi$ of
$G$ (attached to the biholomorphic equivalence $(\sigma,\delta)$) such that
$$
\sigma(g P)=\psi(g)\sigma(P) \quad {\rm for\ all} \ g\in G, P\in C(a,c).
$$
That is, the map $\sigma :C(a,c)\to C(a',c')$ is $G$-twisted-equivariant.

\begin{oss}
We claim that $\psi$ is the identity if  we have a strict covering
isomorphism. The converse does not necessarily hold, 
as it is shown by the example of
$G = \Z / 3\Z$ as a quotient of $ T(3,3,3)$, where all the three 
elements $ \alpha, \beta ,
\gamma$ have the same image $= 1$ mod  $3$ (see the following considerations).
\end{oss}

In order to understand the equivalence relation
induced by the covering isomorphisms on
the set  $\TT(G)$ of triangle structures, recall the following well known facts
from the theory of ramified coverings (see \cite{mi}):

A) The monodromy homomorphism is only determined by the choice of  a
base point  $\infty'$ lying over
$\infty$; a different choice alters the monodromy up to composition with
an inner automorphism (corresponding to a transformation carrying one 
base point to the
other).

B) The map $\delta$ induces isomorphisms
$$\delta_*: \pi_1 (\PP^1_{\C} \setminus B, \infty)\to
\pi_1 (\PP^1_{\C} \setminus B, \delta (\infty))\to
\pi_1 (\PP^1_{\C} \setminus B, \infty)$$
the second being induced by the choice of a path from $\infty$ to
$\delta(\infty)$.

Since the stabiliser of a chosen base point lying over $\infty$ under
the monodromy action is equal to the kernel of the monodromy 
homomorphism $m_\lambda$,
the class of monodromy homomorphisms corresponding to the covering 
$C(a',c')$ is
obtained from the one of the given $\mu$  (corresponding to $C(a,c)$) 
by composing with
$(\delta_*)^{-1}$. In particular, we may set $a'' : = \mu 
(\delta_*)^{-1} (\alpha)$,
and  $c'' : = \mu (\delta_*)^{-1} (\gamma)$. It follows also that 
$\psi$ is gotten from
the natural isomorphism of $\pi_1 (\PP^1_{\C} \setminus B, 
\infty) / ker \mu \ra \pi_1
(\PP^1_{\C} \setminus B,
\infty) / ker (\mu \circ (\delta_*)^{-1})$ induced by $(\delta_*)$, and the
obvious identifications of these quotient groups with $G$ (in more concrete 
terms,
$\psi$ sends $a'' \ra a'$, $c'' \ra c'$).

C) The above shows that if the isomorphism is strict, then $\psi$ is 
the identity.
The converse does not hold  since $\psi$ can be the identity, without
$\delta_*$ being the identity.

\begin{prop}\label{t1}
Let $G$ be a finite group and $(a,c),\, (a',c')\in\TT(G)$. The following
are equivalent:

(i-t) there is a  twisted covering isomorphism from
$\lambda:C(a,c) \rightarrow \PP^1_\C$ to the Galois covering
$\lambda':C(a',c') \rightarrow \PP^1_\C$,

(ii-t) there is a $G$-twisted-equivariant biholomorphic map $\sigma 
:C(a,c)\to C(a',c')$,

(iii-t) $(a,c)$ is in the
${\rm A}_\TT(G)$-orbit of $(a',c')$.

Respectively, the following are equivalent:

(i-s) there is a  strict covering isomorphism from
$\lambda:C(a,c) \rightarrow \PP^1_\C$ to the Galois covering
$\lambda':C(a',c') \rightarrow \PP^1_\C$,

(ii-s) there is a $G$-equivariant biholomorphic map $\sigma 
:C(a,c)\to C(a',c')$,

(iii-s) $(a,c)$ is in the
${\rm I}_\TT(G)$-orbit of $(a',c')$.
\end{prop}

\Proof
The equivalence of (i) and (ii) follows directly from the definition.
In view of A) we only consider triangle structures up to action of 
${\rm Int}(G)$.
We have seen that two triangle structures yield  coverings which are 
twisted covering
isomorphic if and only if there is an automorphism $\delta$ of 
$(\PP^1_{\C} - B)$
and an automorphism $ \psi \in  {\rm Aut} (G)$ such that $\psi \circ \mu  = 
\mu' \circ \delta_*$.

In particular, $a'' : = \mu (\delta_*)^{-1} (\alpha)$,
   $c'' : = \mu (\delta_*)^{-1} (\gamma)$ are $\psi$ equivalent to 
$a', c'$, and it
suffices to show that they are obtained from $ (a,c)$ by one of the 
transformations
$\sigma_i$. Observe however that the group of projectivities
$ {\rm Aut} (\PP^1_{\C} - B)$ is isomorphic to the group of permutations of 
$B$, by the fundamental theorem on projectivities.

We see immediately the action of an element of order two: namely, consider the
projectivity $ z \ra - z$: this leaves the base point $\infty$ fixed, 
as well as the
point $0$, and acts by sending $ \alpha \ra  \gamma, \ \gamma \ra 
\alpha $: we obtain in
this way the transformation $\sigma_3$ on the set of triangle structures.

In order to obtain the transformation $\sigma_1$ of order three, it 
is more convenient,
after a projectivity, to assume that the set $B$ consists of the 
three cubic roots of
unity. Setting $\omega = {\rm exp} ( 2 \pi i /3)$ and  $B = \{ 1, \omega, 
\omega^2\}$,
one sees immediately that $\sigma_1$ is induced by the automorphism
$ z \ra \omega z$, which leaves again the base point $\infty$ fixed 
and cyclically permutes $\alpha, \beta, \gamma$.

\QED

To be able to treat questions of reality we define
\begin{equation}\label{rea1}
\iota(a,c):=(a^{-1},c^{-1})
\end{equation}
for $(a,c)\in \TT(G)$
and call it the {\it conjugate} of $(a,c)$. Note that we have
$\iota(a,c)\in \TT(G)$ and also $\Sigma(\iota(a,c))=\Sigma(a,c)$,
$\mu(\iota(a,c))=\mu(a,c)$.
A feature built into our construction is:

\begin{prop}\label{rea2}
Let $G$ be a finite group and $(a,c)\in \TT(G)$ then
\begin{equation}
C(\iota(a,c))=\overline{C(a,c)}.
\end{equation}
\end{prop}
\Proof For the proof note that by construction the complex conjugates of
the paths $\alpha$, $\gamma$ used in the construction of the triangular curves
$C(a,c)$ satisfy $\bar \alpha=\alpha^{-1}$,  $\bar \gamma=\gamma^{-1}$.
\QED

We now remind the reader of the operations $\sigma_0,\ldots,\sigma_5$
defined in \ref{si1}, \ref{si2}. For later use we observe:

\begin{lem}\label{rea}
Let $G$ be a finite group and $(a,c)\in \TT(G)$.
Let $\rho=\psi\circ\sigma_i\in{\rm A}_\TT(G)$.

(i) In case $i=0$, $\rho(a,c)=\iota(a,c)$ if and only if
$\psi(a)=a^{-1}$ and $\psi(c)=c^{-1}$.

(ii) In case $i=1$, $\rho(a,c)=\iota(a,c)$ if and only if
$\psi(a)=c^{-1}$ and $\psi(c)=ac$.

(iii) In case $i=2$, $\rho(a,c)=\iota(a,c)$ if and only if
$\psi(a)=ac$ and $\psi(c)=a^{-1}$.

(iv) In case $i=3$, $\rho(a,c)=\iota(a,c)$ if and only if
$\psi(a)=c^{-1}$ and $\psi(c)=a^{-1}$.

(v) In case $i=4$, $\rho(a,c)=\iota(a,c)$ if and only if
$\psi(a)=ac$ and $\psi(c)=c^{-1}$.

(vi) In case $i=5$, $\rho(a,c)=\iota(a,c)$ if and only if
$\psi(a)=a^{-1}$ and $\psi(c)=ac$.
\end{lem}

Using the notation of the lemma we assume that
$\rho(a,c)=\psi\circ\sigma_i(a,c)=\iota(a,c)$ and
get the formulae
$$
(\psi\circ\sigma_0)^2(a,c)=(a,c), \
(\psi\circ\sigma_1)^2(a,c)=(c^{-1}ac,c), \
(\psi\circ\sigma_2)^2(a,c)=(a,aca^{-1}),
$$
$$
(\psi\circ\sigma_3)^2(a,c)=(a,c), \
(\psi\circ\sigma_4)^2(a,c)=(c^{-1}ac,c), \
(\psi\circ\sigma_5)^2(a,c)=(a,aca^{-1}),
$$
for the square of $\rho$ on $(a,c)$.

\section{The unmixed case}

In this section we shall translate the
problem of existence and classification of
Beauville surfaces  $S$ of unmixed type into
purely group-theoretic problems.

\subsection{Unmixed Beauville surfaces and group actions}

To have the group theoretic background for the construction of Beauville
surfaces from \cite{cat00} we make the following definition.

\begin{df}
Let $G$ be a finite group $G$. We say that a quadruple $v=(a_1,c_1;a_2,c_2)$
of elements of $G$ is an unmixed Beauville structure for $G$
(for short: u-Beauville)
if and only if

(i) the pairs $a_1,c_1$, and $a_2,c_2$ both generate $G$,

(ii) $\Sigma(a_1,c_1) \cap \Sigma(a_2,c_2) =\{ 1_{G}\}.$

The group $G$ admits an unmixed Beauville structure if such a quadruple $v$
exists. We write $\BBB(G)$ for the set of unmixed Beauville structures on $G$.
\end{df}

We shall need an appropriate notion of equivalence of unmixed
Beauville structures.
In order to clarify it, let us observe that a Beauville surface has
a unique minimal realization (\cite{cat00}, \cite{cat03}), 
and that the Galois group of this covering is isomorphic to $G$.
This yields an action of $G$ on the product $C_1 \times C_2$ (whence, 
two actions of
$G^0$ on both factors) only after we fix an isomorphism of the Galois 
group with $G$.
In turn, these two actions of $G$ determine a triangular covering up to strict
covering isomorphism, so that we can apply Proposition \ref{t1}.

Note that for
$(\psi_1,\psi_2)\in {\rm I}_\TT(G) $
and $(a_1,c_1;a_2,c_2)\in \BBB(G)$ we
have $(\psi_1(a_1,c_1);\psi_2(a_2,c_2))\in \BBB(G)$. Thus we have a faithful
action of ${\rm I}_\TT(G)\times{\rm I}_\TT(G) $ on
$\BBB(G)$. Consider now the group ${\rm B}_\BBB(G)$ generated by the action of
${\rm I}_\TT(G)\times{\rm I}_\TT(G) $ and by the diagonal action of
${\rm Aut}(G) $ (such that
$(\psi\in {\rm Aut}(G) $
carries $(a_1,c_1;a_2,c_2)\in \BBB(G)$ to
$(\psi (a_1,c_1);\psi(a_2,c_2))\in \BBB(G)$).

We additionally define the following
operation
\begin{equation}
\tau((a_1,c_1;a_2,c_2)):=(a_2,c_2;a_1,c_1)
\end{equation}
on the elements of $\BBB(G)$ and let
\begin{equation}
{\rm A}_\BBB(G):=\langle {\rm B}_\BBB(G),\, \tau\rangle
\end{equation}
be the permutation group generated by these permutations of $\BBB(G)$.
Note that ${\rm B}_\BBB(G)$ is a normal subgroup of index $\le 2$ in
${\rm A}_\BBB(G)$.

Given an element $v:=(a_1,c_1;a_2,c_2)\in\BBB(G)$ we define now
\begin{equation}\label{surf}
S(v):=C(a_1,c_1)\times C(a_2,c_2) / G.
\end{equation}
The second condition in our definition of $\BBB(G)$ ensures that the action of
$G$ on the  product has no fixed point, hence the covering
$C(a_1,c_1)\times C(a_2,c_2)\to S(v)$ is unramified.
We call the surface $S(v)$ an {\it unmixed Beauville surface}.
It is obvious
that (\ref{surf}) is a {\it minimal Galois realization} (see \cite{cat00},
\cite{cat03}) of $S(v)$.
Our next result shows that the unmixed Beauville surface
$S(v)$ is isogenous to a {\it higher product}
in the terminology of \cite{cat00}.

\begin{prop}\label{two}
Let $G$ be a finite, non-trivial group with an unmixed Beauville
structure $(a_1,c_1;a_2,c_2)\in\BBB(G)$. Then $\mu(a_1,c_1) < 1$ and
$\mu(a_2,c_2) < 1$. Whence we have: $g(C(a_1,c_1))\ge 2$ and
$g(C(a_2,c_2))\ge 2$.
\end{prop}

\Proof
We may without loss of generality assume that $G$ is not cyclic.
Suppose $(a_1,c_1)$ satisfies $\mu(a_1,c_1)>1$: then the type of $(a_1,c_1)$ is
up to permutation amongst the
$$(2,2,n)\ (n\in \N),\ (2,3,3),\ (2,3,4),\ (2,3,5).$$
In the first case $G$ is a quotient group of the infinite dihedral group and
$G$ cannot admitt an unmixed Beauville structure by Lemma \ref{di}.
There are the following isomorphisms of triangular groups
$$T(2,3,3)=\AAA_4,\ T(2,3,4)=\SSS_4,\ T(2,3,5)=\AAA_5,$$
see \cite{cox}, Chapter 4. These groups do not admit an unmixed Beauville
structure by Proposition \ref{acht}.

If $\mu(a_1,c_1)=1$ then the type of $(a_1,c_1)$ is
up to permutation amongst the
$$(3,3,3),\ (2,4,4),\ (2,3,6)$$
and $G$ is a finite quotient of one of the wall paper groups and cannot admit
an unmixed Beauville structure by the results of Section \ref{wp}.

The second statement follows now from formula \ref{gen} since then
$g(C(a_i,c_i))$, for $i=1,2$, is an integer strictly greater than $1$.
\QED

We may now apply results from
\cite{cat00}, \cite{cat03} to prove:

\begin{prop}\label{eqi}
Let $G$ be a finite group and $v,v'\in\BBB(G)$. Then $S(v)$ is
biholomorphically isomorphic to $S(v')$ if and only if $v$ is in the
${\rm A}_\BBB(G)$-orbit of $v'$.
\end{prop}

\Proof
Let $v=(a_1,c_1;a_2,c_2)$, $v'=(a_1',c_1';a_2',c_2')$.
Assume that there is a biholomorphism between two unmixed Beauville 
surfaces $S(v)$ and
$S(v')$. This happens, by Proposition 3.2
of \cite{cat03}, if and only if there
is a product biholomorphism (up to a possible interchange of the factors)
$$\sigma: C(a_1,c_1)\times C(a_2,c_2) \to C(a_1',c_1')\times C(a_2',c_2')$$
of the product surfaces appearing in
the minimal Galois realization \ref{surf} which normalizes the $G$-action.

In the notation introduced previously, this means that $\sigma$ is twisted
$G$-equivariant. That is, there is an automorphism
$\psi : G\to G$ with $\sigma(g(x,y))=\psi(g)(\sigma(x,y))$ for all $g\in G$
and $(x,y)\in  C(a_1,c_1)\times C(a_2,c_2)$.
Up to replacing one of the two unmixed Beauville structures by an 
${\rm Aut}(G)$ -equivalent
one, we may asume without loss of generality that the map $\sigma$ is strict
$G$-equivariant.

Note that our surfaces are both isogenous to a higher product by Proposition
\ref{two}. Since $\sigma$ is of product type it can
interchange the factors or not. If it does not, there are biholomorphic maps
$$\sigma_1: C(a_1,c_1)\to C(a_1',c_1'), \quad \sigma_2:
C(a_2,c_2)\to C(a_2',c_2')$$
such that $\sigma=(\sigma_1,\sigma_2)$.
If $\sigma$ does interchange the factors there are biholomorphic maps
$$\sigma_1: C(a_1,c_1)\to C(a_2',c_2'),\quad \sigma_2:
C(a_2,c_2)\to C(a_1',c_1')$$
such that $\sigma=(\sigma_1,\sigma_2)$.
In both cases we may now use Proposition \ref{t1} which characterizes strict
$G$-equivariant isomorphisms of triangle coverings.
\QED

\subsection{Unmixed Beauville structures on finite groups}

The question arises: which groups admit Beauville structures?

The unmixed case with $G$ abelian is easy to classify, and all examples were
essentially given in \cite{cat00}, page $24$.

\begin{teo} \label{abelian}
If $G= G^0$  is abelian, non-trivial and admits an
unmixed Beauville structure, then
$ G \cong (\Z/n\Z)^2$, where the integer $n$ is relatively prime to $6$.
Moreover, the structure is critical for both factors.
Conversely, any group $ G \cong (\Z/n\Z)^2$ admits such a structure.
\end{teo}

\Proof Let $(a,c;a',c')$ be an unmixed Beauville structure on $G$, set
$\Sigma:=\Sigma(a,c)$, $\Sigma':=\Sigma(a',c')$ and $b:=a^{-1}c^{-1}$,
$b':={a'}^{-1}{c'}^{-1}$.
Our basic strategy will be to observe that if $H$ is a nontrivial
characteristic
subgroup of $G$, and if we show that for each choice of $\Sigma$
we must have $ \Sigma \supset H $, then we obtain a contradiction to
$ \Sigma \cap \Sigma' = \{ 1\}$.

Consider the primary decomposition of $G$,
$$ G = \bigoplus_{p \in \{{\rm Primes}\}}  G_p $$
and observe that since $G$ is 2-generated, then any $G_p$ (which is a
characteristic subgroup), is also 2-generated.

Step 1.  Let $a = (a_p) \in \bigoplus_{p \in \{{\rm Primes}\}}  G_p$, and let
$\Sigma_p$ be the set of multiples of $a_p, b_p, c_p$: then
$\Sigma \supset \Sigma_p$. This follows since $a_p$ is a multiple of $a$.

Step 2. $G_p \cong (\Z/p^m\Z)^2$.

Since $G_p$ is 2-generated, otherwise $G_p$ is either cyclic
$G_p \cong \Z/p^m\Z$ or $G_p \cong \Z /p^n\Z \oplus \Z/p^m\Z$ with $n < m$.
In both cases the subgroup
$H_p := {\rm Soc} (G_p) := \{ x \in G\, |\, px = 0\}$
is characteristic in $G$ and isomorphic to $\Z /p\Z$.
But $\Sigma \supset \Sigma_p$, and $\Sigma_p$ contains generators of $G_p$,
whence it contains a non-trivial element in the socle, thus
$\Sigma_p \supset H_p := {\rm Soc} (G_p)$, a contradiction.

Step 3. $G_2 = 0$.

Else, by step 2, $G_2 \cong (\Z/ 2^m\Z)^2$, and $H_2 := {\rm Soc} (G_2)
\cong (\Z/ 2\Z)^2$.
But since $\Sigma \supset \Sigma_2$, and $\Sigma_2$ contains a basis
of $G_2$, $\Sigma \supset H_2$, a contradiction.

Step 4. $G_3=0$.

In this case we have that $\Sigma_3$ contains a basis
of $G_3$, whence $\Sigma \cap H_3$ contains at least $6$ nonzero
elements, likewise for $\Sigma' \cap H_3$, a contradiction since
$H_3$ has only $8$ nonzero elements.

Step 5. Whence, $G \cong (\Z/n\Z)^2$, and since $a,b$ are generators of $G$,
they are a basis, and without loss of generality $a,b$ are the standard basis
$e_1, e_2$. It follows that all the elements $a,b,c,a',b',c'$ have order
exactly $n$. Write now the elements of $G$ as row vectors,
$a' := (x,y) , b' := (z,t)$. Then the condition that $\Sigma \cap \Sigma' =
\{1\}$ means that any pair of the six vectors yield a basis of $G$.
By using the primary decomposition, we can read out this condition on
each primary component: thus it suffices to show that there are solutions
in the case where $n = p^m$ is primary.

Step 6. We write up the conditions explicitly, namely, if $n= p^m$
and $ U:= \Z/n\Z^* $,
we want
$$ x,y,z, t \in U, \ x-y,  x+z , z-t, y +t \in U,
x+z -y-t \in U, xt-yz \in U.  $$
Again, these conditions only bear on the residue class modulo $p$,
thus we have $p^{4m-4}$ times the number of solutions that we get
for $n = p$.

Step 7.  Simple counting yields at least $(p-1)(p-2)^2 (p-4)$ solutions.

In this case we get $p-1$ times the number of solutions that we get
for $x=1$, and for each choice of $y \neq 0,1$ $z \neq 0,-1$, $d \neq 0$ we set
$ t := yz + d$ : the other inequalities are then satisfied if $d$ is different
from $ z - yz , -yz - y, (1+z)(1-y)  $  so that the
number of solutions equals at least $(p-1)(p-2)^2 (p-4)$.
\QED

\begin{oss}
The computation above shows that the number of biholomorphism classes
of unmixed Beauville surfaces with abelian group $(\Z/n\Z)^2$ is asymptotic
to at least $ (1/36) \ n^4$ (cf. \cite{bacat} where it is calculated that,
for $n=5$ there are exactly two isomorphism  classes).
\end{oss}

\begin{prop}\label{acht}
No non-abelian group of order $\le 128$ admits an unmixed Beauville
structure.
\end{prop}
\Proof
This result can be obtained by a straightforward computation using 
the computer algebra system MAGMA or by
direct considerations. In fact using the Smallgroups-routine of MAGMA we may
list all groups of order $\le 128$ as explicit permutation-groups or given by
a polycyclic presentation. Loops which are easily designed can be used to 
search for appropriate systems of generators.
\QED

Another simple result is:

\begin{lem}\label{di}
Let G be a non-trivial finite quotient of the infinite dihedral group
D:=$\langle x,y\ | \ x^2,y^2\rangle$
then $G$ does not admit an unmixed Beauville
structure.
\end{lem}
\Proof
The infinite cyclic subgroup $N_0:=\langle xy\rangle$ is normal in 
$D$ of index $2$,
actually $D$ is thus the semidirect product of $N_0 \cong \Z$ through the
subgroup of order $2$ generated by $x$. Let
$t\in D$ be not contained in $N_0$. Then there is an integer $n$ such that
$t=x(xy)^n$.  Since $ y t y = x (xy)^{n-2}$, the normal subgroup 
generated by $t$ then
contains
$(xy)^2$. Hence every normal subgroup $N$ of $D$ not contained in $N_0$
has index $\le 4$ and thus the quotient $D/N$ cannot admit an unmixed Beauville
structure. Let now $N\le N_0$ be a normal subgroup of $D$. The quotient $D/N$
is a finite dihedral group. Let $(a,c)$ be a pair of generators for $D/N$. It
is easy to see that one of the elements $a,\, c,\, ac$ lies in the (cyclic)
image of $N_0$ in $D/N$ and generates it. Thus condition (*) is contradicted.
\QED

\begin{prop}
The following groups admit an unmixed Beauville structure:

1. the alternating groups $\AAA_n$ for large $n$,

2. the symmetric groups $\SSS_n$ for $n\in \N$ with $n\ge 8$ and 
$n\equiv 2$ mod $3$,

3. the groups ${\bf SL}(2,\F_p)$ and  ${\bf PSL}(2,\F_p)$ for every prime
    $p$ distinct from $2,3,5$.
\end{prop}

\Proof 1. Fix two triples $(n_1,n_2,n_3)$, $(m_1,m_2,m_3)\in \N^3$ such that
neither $T(n_1,n_2,n_3)$ nor $T(m_1,m_2,m_3)$ is one of the non-hyperbolic
triangle groups. From \cite{ev} we infer that, for large enough $n\in\N$, the
group $\AAA_n$ has  systems of generators $(a_1,c_1)$ of type 
$(n_1,n_2,n_3)$ and $(a_2,c_2)$ of type $(m_1,m_2,m_3)$. Adding the additional
property that ${\rm gcd}(n_1n_2n_3,m_1m_2m_3)=1$ we find that
$(a_1,c_1;a_2,c_2)$ is an unmixed Beauville structure on $\AAA_n$.

By going through the proofs of \cite{ev} the minimal choice of such an
$n\in\N$ can be made effective.

2. This follows directly from the first Proposition of Section 5.1.

3. Let $p$ be a prime with the property that
   no prime $q\ge 5$ divides $p^2-1$. Then $(p,1)$ is a primitive solution 
of an equation  
\begin{equation}\label{sie}
y^2-x^3=\pm 2^n3^m
\end{equation}
with $n,m\in \N$ chosen appropriately. It is known that the collection of these
equations has $98$ primitive solutions (as $n,m$ vary). A table of them is
contained in \cite{bir} Table 4, page 125. From this we see that
$p=2,3,5,7,17$ are the only primes with the property that
no prime $q\ge 5$ divides $p^2-1$. Notice that a theorem of C.L. Siegel 
implies directly that there are only finitely many such primes.
A special case of this theorem says that any of the equation (\ref{sie}) has
only finitely many solutions in $\Z[1/2,1/3]$. 

If $p$ is a prime with $p\ne 2,3,5,7,17$
we use the system of
   generators from (\ref{ps1}) which is of type $(4,6,p)$ together with one of
   the system of generators from (\ref{ps2}) or (\ref{ps3}) 
to conclude the result
   for ${\bf SL}(2,\F_p)$. The groups  ${\bf PSL}(2,\F_p)$ can be treated by
   reduction of these systems of generators, observing that the two
generators belonging to  different  systems have coprime orders.

For the primes $p=7,17$ appropriate systems of generators can be easily found
by a computer calculation.

\QED

From the third item of the above proposition we immediately obtain a 
{\it proof of Theorem \ref{intro1}}.

\begin{oss}
The various systems of generators given in Section 5 for the alternating 
groups $\AAA_n$ and for ${\bf SL}(2,\F_p)$, ${\bf PSL}(2,\F_p)$ can be grouped
together in many ways to construct unmixed Beauville structures on these
groups. We in turn obtain Beauville surfaces of unmixed type for which the two
curves appearing in  the minimal Galois realization have different genus.
\end{oss}

We shall describe now some groups of a completely different nature
admitting an unmixed Beauville structure.
For $n\in \N$ define:
\begin{equation}\label{x}
C[n]:=\langle x,y \ | \
x^n=y^n=(xy)^n=(xy^{-1})^n=(xy^{-2})^n=(xy^{-1}xy^{-2})^n=1\rangle 
\end{equation}

\begin{prop}
Let $n\in \N$ be given with $ G.C.D. (n,6)=1$. Let further $N\le C[n]'$ be a 
normal subgroup of finite index where $C[n]'$ is the commutator subgroup of 
$C[n]$. Then $C[n]/N$ admits an unmixed
Beauville structure.
\end{prop}

\Proof
An unmixed Beauville structure for $C[n]/N$ is given by
$$(a_1,c_1;a_2,c_2)=(x,y;xy^{-1},xy^{-2}).$$ 

In fact, let $G^{ab}$ be the abelianization of $G$. Then $\Sigma$ and $\Sigma'$
map injectively into $G^{ab}$, and their images do not meet inside $G^{ab}$,
as verified in \cite{cat00}, lemma 3.21.

\QED

 Amongst the quotients $C[n]/N$ ($N\le
C[n]'$) are all finite groups $G$ of exponent $n$ having 
$\Z/n\Z\times \Z/n\Z$ as abelianization. The Proposition can hence be
used to construct finite $p$-groups ($p\ge 5$) admitting an unmixed
Beauville structure.

\subsection{Questions of reality}

We shall also translate into group theoretic conditions the
two questions concerning an unmixed Beauville surface
mentioned in the introduction:

\begin{itemize}
\item
Is $S$  biholomorphic to the complex conjugate surface $\bar{S}$?
\item
Is $S$ real, i.e. does there exist such a biholomorphism $\sigma$ with the
property that $\sigma^2 = {\rm Id}$?
\end{itemize}

Let $G$ be a finite group and  $v=(a_1,c_1;a_2,c_2)\in\BBB(G).$
In analogy with  (\ref{rea1}) we define
\begin{equation}
\iota(v):=(a_1^{-1},c_1^{-1};a_2^{-1},c_2^{-1})
\end{equation}
and infer from Proposition \ref{rea2}:
\begin{equation}
S(\iota(v))=\overline{S(v)}.
\end{equation}
 From Proposition \ref{eqi} we get

\begin{prop}\label{rea12}
Let $G$ be a finite group with an unmixed Beauville structure
$v=(a_1,c_1;a_2,c_2)\in\BBB(G).$ Then

1. $S(v)$  biholomorphic to $\overline{S(v)}$ if and only if $\iota(v)$ is
in the ${\rm A}_\BBB(G)$ orbit of $v$,

2. $S(v)$ is real if and only if there is a $\rho\in {\rm A}_\BBB(G)$ with
$\rho(v)=\iota(v)$ and moreover $\rho(\iota(v))=v$.
\end{prop}

\begin{oss}
The above observations immediately imply that unmixed Beauville surfaces
$S$ with abelian group $G$ always have a real structure, since $g \ra -g$
is an automorphism (of order $2$).
\end{oss}

We observe the following:
\begin{cor}\label{rea13}
Let $G$ be a finite group with an unmixed Beauville structure
$v=(a_1,c_1;a_2,c_2)\in\BBB(G).$
Assume that the sets $\{ {\rm ord}(a_1), {\rm ord}(c_1),{\rm ord}(a_1c_1)\}$
and $\{ {\rm ord}(a_2), {\rm ord}(c_2),{\rm ord}(a_2c_2)\}$ are distinct.
Assume further that both $(a_1,c_1)$ and $(a_2,c_2)$ are strict.
Then $S(v)$ is biholomorphically isomorphic to $\overline{S(v)}$
if and only if the following holds:

There are  inner automorphisms $\phi_1,\phi_2$ of
$G$ and an automorphism $\psi \in Aut(G)$
such that, setting
$\psi_j : = \psi \circ \phi_j$, we have
$\psi_1(a_1) = a_1^{-1}$, $\psi_1 (c_1) = c_1^{-1}$,
   and $\psi_2 (a_2) =
   {a_2}^{-1}$, $\psi_2 (c_2) = {c_2}^{-1}$.

In particular  $S(v)$ is isomorphic to $\overline{S(v)}$ if and only if
$S(v)$ has a real structure.
\end{cor}

\Proof The first statement follows from our Definition
of ${\rm A}_\BBB(G)$ and Proposition \ref{eqi}. In fact let $\rho\in
{\rm A}_\BBB(G)$ be such that $\rho(v)=\iota(v)$. We have
$$\rho=(\psi_1\circ \sigma_i,\psi_2\circ \sigma_j)\circ \tau^e$$
with $\psi_1,\, \psi_2$  as above, $i,j\in\{0,\ldots,5\}$ and
$e\in\{0,1\}$.
Our incompatibility conditions on the orders imply that  $e=0$
and $i=j=0$ (see Lemma \ref{rea}).

For the second statement note that the conclusion implies that
both $\psi_1$ and $\psi_2$ have order $2$. \QED

\begin{oss}
If the unmixed Beauville structure $v$ does not have the
strong incompatibility properties of the corollary then
lemma \ref{rea} gives the appropriate conditions.
\end{oss}

 From our Corollary we immediately get:

{\it Proof of Theorem 1.2.}
Let $v:=(a,c;a',c')$ with $a,c,a',c'\in\SSS_n$ as in Proposition
\ref{sy11}. Then $v$ is an unmixed Beauville structure on $\SSS_n$.
The type of $(a,c)$ is $(6, 3(n-3),3(n-4))$, while the type of $(a',c')$ is
$(n, n-1, n)$ or $(n, n-1, (n^2-1)/4)$.
Suppose that $S(v)$ is biholomorphic to $\overline{S(v)}$.
By Proposition \ref{rea12} (statement 1) $\iota(v)$ is
in the ${\rm A}_\BBB(\SSS_n)$ orbit of $S(v)$. The incompatibility of 
the types of
$(a,c)$ and $(a',c')$ makes Corollary  \ref{rea13} applicable. This 
implies that there
is a $\psi\in {\rm Aut}(\SSS_n)$ with $\psi(a)=a^{-1}$ and
$\psi(c)=c^{-1}$. Since all automorphisms on $\SSS_n$ ($n\ge 8$) are inner
we obtain a contradiction to Proposition \ref{sy11}.
\QED

As noted above the unmixed Beauville surfaces
$S$ coming from an abelian group $G$ always have a real structure. It 
is also possible to construct examples from non-abelian groups:

\begin{prop} Let $p\ge 5$ be a prime with $p\equiv 1$ mod $4$. Set $n:=3p+1$.
Then there is an unmixed Beauville structure $v$ for the group $\AAA_n$ 
such that
$S(v)$ is biholomorphic to $\overline{S(v)}$.
\end{prop}
\Proof We use the first and the second system of generators for $\AAA_n$
from Proposition \ref{alp2}. The first is a system of generators 
$(a_1,c_1)$ of type
$(2,3,84)$, the second gives $(a_2,c_2)$ of type $(p,5p,2p+3)$.
Since the orders in the two types are coprime $(a_1,c_1;a_2,c_2)$ is an unmixed
Beauville structure on $\AAA_n$.
The existence of the respective elements $\gamma$ in Proposition \ref{alp2}
implies the last assertion. Note that both the elements $\gamma$ can be chosen
to be in $\SSS_n\setminus \AAA_n$. \QED

In the further arguments we shall often use the fact that every automorphism
of $\AAA_n$ ( $n\ne 6$) is induced by 
conjugation by an element of $\SSS_n$ (see \cite{suz}, page 299).

\bigskip

\begin{prop}
The following groups admit unmixed Beauville structures $v$ such that $S(v)$ is
not biholomorpically isomorphic to $\overline{S(v)}$:

1. the symmetric group $\SSS_n$ for $n\ge 8$ and $n\equiv 2$ mod $3$,

2. the alternating group $\AAA_n$ for $n\ge 16$ and $n\equiv 0$ mod $4$,  
$n\equiv 1$ mod $3$, $n\not\equiv 3,4$ mod $7$.

\end{prop}

\Proof 
1. This is just the example of Section 5.1.

2. We use the system of generators $(a_1,c_1)$ from Proposition
\ref{alp2}, 1. It has type $(2,3,84)$. We then choose $p=5$ and $q=11$ and get
from Proposition \ref{alp1} a system of generators $(a_2,c_2)$ of type
$(11,5(n-11),n-3)$. Both systems are strict. We set $v:=(a_1,c_1;a_2,c_2)$.
The congruence conditions $n\not\equiv 3,4$ mod $7$ insure that 
$\nu(a_1,,c_1)$ is coprime to $\nu(a_2,,c_2)$, hence
this is an
unmixed Beauville structure. It also satisfies the hypotheses of Corollary
\ref{rea13}. If $S(v)$ is
biholomorpically isomorphic to $\overline{S(v)}$ we obtain an element
$\gamma \in\SSS_n$ with $\gamma a_2\gamma^{-1}=a_2^{-1}$ and  
 $\gamma c_2\gamma^{-1}=c_2^{-1}$. This contradicts Proposition \ref{alp1}.
\QED

\begin{prop}
Let $p>5$ be a prime with $p\equiv 1$ mod $4$, $p\not\equiv 2,4$ mod 5,  
$p\not\equiv 5$ mod $13$ and $p\not\equiv 4$ mod $11$. Set $n:=3p+1$. Then 
the alternating group $G: = \AAA_n$ admits an unmixed Beauville structure $v$
such that there is an 
element $\alpha\in {\rm A}_\BBB(G)$ with 
$\alpha(v)=\iota(v)$ but such that there is no element 
$\beta\in {\rm A}_\BBB(G)$ with $\beta(v)=\iota(v)$ and 
$\beta(\iota(v))=v$.
\end{prop}

\Proof In order to construct $v$ we use the system of 
generators $(a_1,c_1)$ from Proposition
\ref{alp3}. It has type $(3p-2,3p-1,3p-1)$. We then use the system of
generators $(a_2,c_2)$ from Proposition \ref{alp2}, 2. It has type 
$(p,5p,2p+3)$. The second system is strict. We set $v:=(a_1,c_1;a_2,c_2)$.
The congruence conditions $p\not\equiv 2,4$ mod 5,  
$p\not\equiv 5$ mod $13$ and $p\not\equiv 4$ mod $11$ ensure that
$\nu(a_1,,c_1)$ is coprime to $\nu(a_2,,c_2)$, hence $v$ is an
unmixed Beauville structure.

We shall first show that there exists $\alpha\in {\rm A}_\BBB(\AAA_n)$ with 
$\alpha(v)=\iota(v)=(a_1^{-1},c_1^{-1};a_2^{-1},c_2^{-1})$. We choose 
$\gamma_1$ as in Proposition \ref{alp3} and 
$\gamma_2$ as in Proposition \ref{alp2}, 2. Let $w\in\SSS_n$ be a
representative of the nontrivial coset of $\AAA_n$ in $\SSS_n$.
By Propositions \ref{alp3},  \ref{alp2} these choices can be made so that
$\gamma_1=\delta_1w$, $\gamma_2=\delta_2w$ with $\delta_1,\delta_2\in \AAA_n$.
We have now 
\begin{equation}
(w^{-1}\delta_1^{-1} \ a_1 \ \delta_1w, 
w^{-1}\delta_1^{-1}\  c_1^{-1}a_1^{-1} \ \delta_1w)=(a_1^{-1},c_1^{-1}),
\end{equation}
\begin{equation}
(w^{-1}\delta_2^{-1}\  a_2\  \delta_2 w, 
w^{-1}\delta_2^{-1}\  c_2\  \delta_2w)=
(a_2^{-1},c_2^{-1}).
\end{equation}
Recalling the formula for $\sigma_5$ (see (\ref{si2})) the 
existence of $\rho$ follows from
our definition of
${\rm A}_\BBB(\AAA_n)$.

Suppose now that there is a $\beta\in {\rm A}_\BBB(G)$ as indicated, then
$\beta^2(v)=v$. By construction of $v$ the transformation $\beta$ cannot
interchange $(a_1,c_1)$ and $(a_2,c_2)$. Hence we find $\beta_1,\beta_2\in 
{\rm A}_\TT(\AAA_n)$ with 
$$\beta_1(a_1,c_1)=(a_1^{-1},c_1^{-1}),\quad \beta_1(a_2,c_2)=
(a_2^{-1},c_2^{-1}).$$ 
from $\beta_1^2(a_1,c_1)=(a_1,c_1)$ and the formulae given immediately after
Lemma \ref{rea} we infer that either $\beta_1=\psi\circ\sigma_0$ or 
$\beta_1=\psi\circ\sigma_3$ for a suitable automorphism $\psi$ of $\AAA_n$. 
(Note that $a_1$ and $c_1$ cannot commute.)
Going back to Lemma \ref{rea} (i), (iv) we find a contradiction against the
statement of Proposition \ref{alp3}.\QED
 
{\it Theorem \ref{intro3}} follows immediately from the above Proposition and
from Proposition \ref{rea12}.

\section{The mixed case}

In this section we will first fix the algebraic data that
are needed for the construction of  Beauville surfaces of mixed type.
Later on we will use this description to give several examples.

\subsection{Mixed Beauville surfaces and group actions}

This subsection contains the translation between the geometrical data 
of a mixed
Beauville  surface and the corresponding algebraic data:  finite 
groups endowed with a
mixed  Beauville structure. This concept is contained in the following:

\begin{df}

Let $G$ be a non-trivial finite group. A mixed Beauville
quadruple for $G$ is a quadruple
$M=(G^0;a,c;g)$ consisting of a
subgroup $G^0$ of index $2$ in $G$, of elements
$a,c\in G^0$ and of an element $g\in G$ such that

1. $G^0$ is generated by $a,c$,

2. $g\notin G^0$,

3. for  every $\gamma \in G^0$ we have
$g\gamma g\gamma\notin \Sigma(a,c)$,

4. $\Sigma(a,c)\cap\Sigma(gag^{-1},gcg^{-1})=\{ 1_G \}$.

 From a mixed Beauville quadruple we obtain, by forgetting about the 
choice of $g$, a
mixed Beauville  triple for $G$, $u =(G^0;a,c)$.
The group $G$ is said to admit a mixed Beauville structure if such a 
quadruple $M$
exists. We let then $\M_4(G)$ be the set of mixed Beauville quadruples on the
group $G$,  $\M_3(G)$ be the set of mixed Beauville triples on the
group $G$. These last will also be called mixed Beauville structures.
\end{df}

We shall describe now the correspondence between the data for an 
unmixed Beauville
structure given above and those given in \cite{cat00}
(also described in the introduction).
Let $M=(G^0;a,c;g)$ be a mixed Beauville quadruple on a finite group $G$.
Then $G^0$ is normal in $G$. By condition (3) the exact sequence
\begin{equation}\label{exa}
1 \rightarrow G^0 \rightarrow G \rightarrow \mathbb{Z} / 2 \Z \rightarrow 1,
\end{equation}
does not split. Define
$\varphi_g :G^0\to G^0$ to be
the automorphism of $G^0$ induced by conjugation with
$g$, that is $\varphi_g(\gamma)=g\gamma g^{-1}$ for all $\gamma\in G^0$.
Suppose $\varphi_g$ would be an inner automorphism. Then we could find
$\delta \in G^0$ with $\varphi_g(\gamma)=\delta\gamma\delta^{-1}$ for all
$\gamma\in G^0$. This implies that $\Sigma(gag^{-1},gcg^{-1})=\Sigma(a,c)$.
Since $G$ is nontrivial condition (4) cannot hold.

Let $ \tau := \tau_g := g^2 \in G^0$. We
have $\varphi_g (\tau ) = \tau$ and $\varphi_g^2 ={\rm Int}_{\tau}$ where
${\rm Int}_{\tau}$ is the inner automorphism induced by $\tau$.
This shows that
$\varphi_g$ is of order $2$ in the group of outer automorphisms
${\rm Out}(G^0)$ of $G^0$. Conversely given a non-trivial finite group $G^0$
together with an an automorphism $\varphi :G^0\to G^0$ of order $2$ in the
outer automorphism group allows us to find a group $G$ together with an exact
sequence (\ref{exa}).

It is important to observe that the conditions (3), (4) are the ones 
which guarantee
the freeness of the action of $G$.

We shall describe now the appropriate notion of equivalence for mixed
Beauville structures. Let $M=(G^0;a,c;g)$ be a mixed Beauville quadruple for
the group $G$ and $\psi:G\to G$ be an automorphism of $G$: then
$\psi(M):=(\psi(G^0);\psi(a),\psi(c);\psi(g))$
is again a mixed Beauville structure
on $G$. Thus we obtain respective actions of ${\rm Aut}(G)$ on  $\M_4(G)$,
$\M_3(G)$.
If $M=(G^0;a,c;g)$ is a mixed Beauville quadruple for
the group $G$ and $\gamma\in G^0$ then $M_\gamma=(G^0;a,c;\gamma g)$ is also
  a mixed Beauville quadruple on $G$.

We can therefore, without loss of generality, only consider   mixed Beauville
triples (beware, such a triple is obtained from a quadruple 
satisfying  conditions
(1)-(4) of the previous definition).

We consider  on the set $\M(G) : = \M_3(G)$ of mixed Beauville structures
the action of the group
\begin{equation}
{\rm A}_{\M}(G):=< {\rm Aut}(G),\
\sigma_3,\sigma_4>.
\end{equation}
with the understanding that the operations $\sigma_3,\sigma_4$ from
(\ref{si1}), (\ref{si2}) are applied to the pair $(a,c)$ of 
generators of $G^0$.
Note that the operations $\sigma_1,\sigma_2,\sigma_5$ are also in
${\rm A}_{\M}(G)$ because of (\ref{si4}).

We will recall now how the above algebraic data give rise to a Beauville
surface of mixed type. Let $u:=(G^0;a,c;g)$ be a mixed Beauville quadruple on
$G$. Set $\tau_g:=g^2$ and
$\varphi_g(\gamma):= g\gamma g^{-1}$ for $\gamma\in G^0$.
By Riemann's existence theorem as in the previous section
the elements $a$, $b=a^{-1}c^{-1}$, $c$ give rise to a
Galois covering $\lambda: C(a,c) \rightarrow \mathbb{P}^1_\C$ ramified only in
$\{-1, 0, 1 \}$ with ramification indices equal to the respective
orders of $a$, $b=a^{-1}c^{-1}$, $c$ and with group
$G^0$. We let the group $G$ act on
$C(a,c) \times C(a,c)$ by letting
\begin{equation}
\gamma (x,y) = (\gamma x, \varphi_g(\gamma) y), \ \ \ g (x,y) = (y, \tau_g x),
\end{equation}
for all $\gamma \in G^0$ and
$(x,y) \in C(a,c) \times C(a,c)$. These  formulae determine an action
of $G$ uniquely.
By our conditions (3), (4) in the definition of a mixed Beauville quadruple on
$G$ the above action of $G$ is fixed point free, yielding a Beauville surface
of mixed type
\begin{equation}\label{surfm}
S(u) := C(a,c) \times C(a,c) /G.
\end{equation}

It is obvious
that (\ref{surfm}) is a {\it minimal Galois representation} (see \cite{cat00},
\cite{cat03}) of $S(u)$. From Proposition \ref{two} we infer
that the mixed Beauville surface
$S(u)$ is isogenous to a {\it higher product}
in the terminology of \cite{cat03}.

Observe that a Beauville surface of mixed type
$S(u) = C(a,c) \times C(a,c)/G$
has a natural unramified
double cover $S^0(u) = (C(a,c) \times C(a,c))/G^0$ which is of unmixed type.

\begin{prop}\label{eqi5}
Let $G$ be a finite group and $u_1,u_2\in\M(G)$. Then $S(u_1)$ is
biholomorphic to $S(u_2)$ if and only if $u_1$ is in the
${\rm A}_{\M}(G)$-orbit of $u_2$.
\end{prop}
\Proof This Proposition follows (as Proposition \ref{eqi}) from \cite{cat00},
\cite{cat03} (Proposition 3.2).

Let $M=(G^0_1,a_1,c_1;g_1)$, $M'=(G^0_2;a_2,c_2;g_2)$.
Assume that the two unmixed Beauville surfaces $S(u)$ and $S(u')$ are
biholomorphically isomorphic. This happens, by Proposition 3.2
of \cite{cat03}, if and only if there
is a product biholomorphism (up to an interchange of the factors)
$$\sigma: C(a_1,c_1)\times C(a_1,c_1) \to C(a_2,c_2)\times C(a_2,c_2)$$
of the product surfaces.
Since $\sigma$ is of product type it can
interchange the factors or not. Hence there are biholomorphic maps
$$\sigma_1,\, \sigma_2: C(a_1,c_1)\to C(a_2,c_2) $$
with
$$\sigma(x,y)=(\sigma_1(x),\sigma_2(y))\qquad {\rm for\ all}\
x,y\in C(a_1,c_1)$$
in case $\sigma$ does not interchange the factors and
$$\sigma(x,y)=(\sigma_1(y),\sigma_2(x))\qquad {\rm for\ all}\
x,y\in C(a_1,c_1)$$
in case $\sigma$ does interchange the factors.
The map $\sigma$ normalises the $G$ action if there is an automorphism
$\psi : G\to G$ with $\sigma(g(x,y))=\psi(g)(\sigma(x,y))$ for all $g\in G$
and $(x,y)\in  C(a_1,c_1)\times C(a_1,c_1)$.

In both cases we may now use Proposition \ref{t1} together
with some straightforward computations to complete the only if
statement of our Proposition.

The reverse statement follows from Proposition
\ref{t1} together with the apropriate considerations.
\QED

\subsection{Mixed Beauville structures on finite groups}

To find a group $G$ with a mixed Beauville structure is rather difficult,
for instance the subgroup $G^0$ cannot be abelian:

\begin{teo} \label{nonabelian}
If a group $G$ admits a mixed Beauville structure, then the subgroup
$G^0$ is non abelian.
\end{teo}

\Proof
By  theorem \ref{abelian} we know that $G^0$ is isomorphic to
$(\Z/n\Z)^2$, where $n$ is an odd number not divisible by $3$.

In particular, multiplication by $2$ is an isomorphism of $G^0$, thus
there is a unique element $\gamma$ such that $ - 2 \gamma = \tau$.
Since $ - 2 \varphi (\gamma) = \varphi (\tau) = \tau = - 2 \gamma$, it follows
that $  \varphi (\gamma) = \gamma $, and we have found a solution to the
prohibited equation $ \varphi (\gamma) + \tau + \gamma \in \Sigma$,
since $ 0 \in \Sigma$. Whence the desired contradiction.
\QED

We also report the following fact obtained by computer calculations using
MAGMA.

\begin{prop}
No group of order $\le 512$ admits a mixed Beauville structure.
\end{prop}

We shall describe now a general construction which gives finite groups $G$
with a mixed Beauville structure. Let $H$ be non-trivial group. Let
$\Theta :H\times H \to H\times H$ be the automorphism defined by
$\Theta(g,h):=(h,g)$ ($g,h\in H$). We consider the semidirect product
\begin{equation}
H_{[4]}:= (H\times H) \rtimes  \Z / 4\Z
\end{equation}
where the generator $1$ of $ \Z / 4\Z$ acts through $\Theta$ on $H\times H$.
Since $\Theta^2$ is the identity we find
\begin{equation}
H_{[2]}:=H\times H\times 2\Z / 4\Z\cong H\times H\times \Z / 2\Z
\end{equation}
as a subgroup of index $2$ in $H_{[4]}$.

Notice that the exact sequence
$$
1 \rightarrow H_{[2]} \rightarrow H_{[4]} \rightarrow
\mathbb{Z} / 2 \Z \rightarrow 1,
$$
does not split because there is no element of order $2$ in $ H_{[4]}$ which
is not already contained in $ H_{[2]}$.

We have
\begin{lem}\label{vz3}
Let $H$ be a non-trivial group and let $a_1,c_1$,  $a_2,c_2$  be elements of
$H$. Assume that

1. the orders of $a_1,c_1$ are even,

2. $a_1^2, a_1c_1, c_1^2$ generate $H$,

3. $a_2,c_2$ also generate $H$,

4. $\nu(a_1,c_1)$ is coprime to  $\nu(a_2,c_2)$.

Set $G:=H_{[4]}$, $G^0:=H_{[2]}$ as above
and $a:=(a_1,a_2,2)$, $c:=(c_1,c_2,2)$. Then
$(G^0;a,c)$ is a mixed Beauville structure on $G$.

If $G$ is a perfect group then the conclusion holds with hypothesis 2 replaced
by the following hypothesis

2'. $a_1,c_1$ generate $H$.
\end{lem}

\Proof We first show that $a,c$ generate  $G^0:=H_{[2]}$.
Let $L:=\langle a,c\rangle$. We view $H\times H$
as the subgroup  $H\times H \times \{ 0\}$ of $H_{[2]}$. The elements
$a^2,ac,c^2$ are in this subgroup. Condition 2 implies that $L\cap (H\times H)$
projects surjectively onto the first factor of $H\times H$. From conditions
1,3, 4 we infer that $a_2, c_2$ have odd order, and that there is an 
even number
$2m$ such that
  $a_2^{2m}, c_2^{2m}$  generate $H$, while $a_1^{2m}=  c_1^{2m} = 1$. 
It follows that
$ H \times H\le L$, and it is then obvious that
$ L =G^0$.

Observe next that
\begin{equation}\label{crit}
(1_H,1_H,2)\notin \Sigma(a,c).
\end{equation}
It would have to be  conjugate of a power of $a$ , $c$ or $b$. Since the
orders of $a_1$, $b_1$, $c_1$ are even, we obtain a contradiction. Note in
fact that the third component of $ac$ is $0$ by construction. 

We shall now verify the third condition of our definition of a
mixed Beauville structure. Suppose preliminarly that 
$h=(x,y,z)\in\Sigma(a,c)$ satisfies
${\rm ord}(x)={\rm ord}(y)$: then our condition 4 implies that $x=y=1_H$ and
(\ref{crit}) shows $h=1_{H_{[4]}}$.

Let now $g\in H_{[4]}$, $g\notin H_{[2]}$
and $\gamma\in G^0=H_{[2]}$ be given. Then $g\gamma=(x,y,\pm 1)$ for
appropriate $x,y\in H$. We find
$$ (g\gamma)^2=(xy,yx,2)$$
and the orders of the first two components of $(g\gamma)^2$ are the same.
The remark above shows that $(g\gamma)^2\in\Sigma(a,c)$ implies 
$(g\gamma)^2=1$.

We come now to the fourth condition of our definition of a
mixed Beauville quadruple. Let $g\in H_{[4]}$, $g\notin H_{[2]}$ be given,
for instance $(1_H,1_H,1)$.
Conjugation with $g$ interchanges then
the first two components of an element $h\in H_{[4]}$.
Our hypothesis 4 implies the result.

So far we have proved the lemma using hypothesis 2. Assume that $H$ is a 
perfect group (this means that $H$ is generated by commutators). Because of
hypothesis 2' the group $H$ is generated by commutators of words in 
$a_1,c_1$.  Defining $L$ as before we see again that that $L\cap (H\times H)$
projects surjectively onto the first factor of $H\times H$. The rest of the
proof is the same. 

\QED

As an application we get

\begin{prop}
Let $H$ be one of the following groups:

1. the alternating group $\AAA_n$ for large $n$,

2. ${\bf SL}(2,\F_p)$ for $p\ne 2,3,5,17$.

Then $H_{[4]}$ admits a mixed Beauville structure.

\end{prop}

\Proof 1. Fix two triples $(n_1,n_2,n_3)$, $(m_1,m_2,m_3)\in \N^3$ such that
neither $T(n_1,n_2,n_3)$ nor $T(m_1,m_2,m_3)$ is one of the non-hyperbolic
triangle groups. From \cite{ev} we infer that, for large enough $n\in\N$, the
group $\AAA_n$ has  systems of generators $(a_1,c_1)$ of type 
$(n_1,n_2,n_3)$ and $(a_2,c_2)$ of type $(m_1,m_2,m_3)$. Adding the additional
properties that $n_1,n_2$ are even and 
${\rm gcd}(n_1n_2n_3,m_1m_2m_3)=1$ we find that the 
$(a_1,c_1;a_2,c_2)$ satisfy the hypotheses 1,2',3,4 of the previous
lemma. Since $\AAA_n$ is, for large $n$, a simple group the statement follows.

2. The primes $p=2,3,5,17$ are the only primes with the property that
   no prime $q\ge 5$ divides $p^2-1$. In the other cases we use the system of
   generators from \ref{ps1} which is of type $(4,6,p)$ together with one of
   the system of generators from \ref{ps2} or \ref{ps3} to obtain 
   generators satisfying  hypotheses 1,2',3,4 of the previous lemma. Since
   ${\bf SL}(2,\F_p)$ is a perfect group (for $p\ne 2,3$) the statement
   follows. 

\QED

\subsection{Questions of reality}

Let $G$ be a finite group and  $u=(G^0;a,c)\in \M(G) = \M_3(G).$
In analogy with \ref{rea1} we define
\begin{equation}
\iota(u):=(G^0;a^{-1},c^{-1})
\end{equation}
and infer from Proposition \ref{rea2}:
\begin{equation}
S(\iota(u))=\overline{S(u)}.
\end{equation}
 From Proposition \ref{eqi} we get

\begin{prop}\label{rea14}
Let $G$ be a finite group and $u\in\M(G)$ then

1. $S(u)$  biholomorphic to $\overline{S(u)}$ if and only if $\iota(u)$ lies
in the ${\rm A}_{\M}(G)$ orbit of $u$,

2. $S(u)$ real if and only if there exists $\rho\in {\rm A}_{\M}(G)$ with
$\rho(u)=\iota(u)$ and $\rho(\iota(u))=u$.
\end{prop}

Observe that if a mixed Beauville surface $S$ is isomorphic to its conjugate,
then necessarily the same holds for its natural unmixed double cover $S^0$.

We will now formulate an algebraic condition on $u\in\M(G)$ which will allow
us to show easily
that the associated Beauville surface $S(u)$ is not isomorphic
to $\overline{S(u)}$.

\begin{cor}\label{vd1}
Let $G$ be a finite group and  $u=(G^0;a,c)\in\M(G)$
and assume that $(a,c)$ is a strict system of generators for $G^0$.
Then $S(u) \cong \overline{S(u)}$ if and only if
there is an automorphism $\psi$ of $G$ such that $\psi(G^0)=G^0$ and
   $\psi (a) = a^{-1}$, $\psi (c) = c^{-1}$.
\end{cor}

\Proof
This Corollary follows from Proposition \ref{rea14} in the same way as
Corollary \ref{rea13} follows from Proposition \ref{rea12}.

\QED

We shall now give a alternative description of the conclusion of Corollary
\ref{vd1}. 

\begin{oss} With the assumptions of Corollary \ref{vd1} let $g\in G$ represent
  the nontrivial coset of $G^0$ in $G$. Set $\tau:=\tau_g=g^2$ and $\varphi:=
\varphi_g$. Then $S(u) \cong \overline{S(u)}$ if and only
if  there is an automorphism $\beta$ of $G^0$ such that
$\beta (a) = a^{-1}$, $\beta (c) = c^{-1}$, and an element $\gamma \in G^0$
such that $ \tau (\beta (\tau^{-1} ) ) = \varphi (\gamma ) \gamma$.
\end{oss}

\Proof
$S \cong \bar{S}$ if and only if $C_1 \times C_2$ admits an antiholomorphism
$\sigma$ which normalizes the action of $G$.
Since there are biholomorphisms of $C_1 \times C_2$ which exchanges the factors
(and lies in $G$), we may assume that such an antiholomorphism does not
exchange the two factors. Being of product type $\sigma = \sigma_1 \times
\sigma_2$, it must normalize the
product group $ G^0 \times G^0$. We get thus a pair of automorphisms 
$\beta_1, \beta_2$ of $G$.
Since $\beta_1 \times \beta_2$ leaves the subgroup
$ \{ (\gamma, \varphi(\gamma)) \ | \ \gamma\in G \}$ invariant , 
it follows that
$ \beta_2 = \varphi \beta_1 \varphi^{-1} $, and in particular $\beta_2$ 
carries
$a' := \varphi (a), c' := \varphi (c)$ to their respective inverses.

Now, $\sigma_1 \times \sigma_2$  normalizes the whole subgroup $G$
if and only if for each $\epsilon \in G^0$ there is $\delta \in G^0$ such that
$$ \sigma_1 \, \varphi (\epsilon) \sigma_2^{-1}= \varphi (\delta) 
   \sigma_2 (\tau \epsilon) \sigma_1^{-1}= \tau \delta .$$
We use now the strictness of the structure: this ensures that both
$\sigma_i$'s are liftings of the standard complex conjugation, whence
we easily conclude that there is an element $ \gamma \in G^0$ such that
$\sigma_2 = \gamma \sigma_1$.

From the second equation we conclude that
$ \delta  = \tau^{-1} \gamma \sigma_1 \tau \epsilon \sigma_1^{-1}$, and
the first boils then down to $ \sigma_1 (\varphi (\epsilon)) \sigma_1^{-1}
{\gamma} ^{-1} =  \tau^{-1} (\varphi(\gamma) ) \gamma \sigma_1 
\tau (\varphi (\epsilon))
\sigma_1^{-1} {\gamma}^{-1}$.

Since this must hold for all $\epsilon \in G^0$, it is equivalent to require
$ \sigma_1 \tau^{-1} \sigma_1^{-1} = \tau^{-1} (\varphi (\gamma) ) \gamma$,
i.e., $ \tau (\beta (\tau^{-1} ) ) = (\varphi (\gamma )) \gamma$.

\QED

\bigskip

We shall now give examples of mixed Beauville structures. In the proofs we
shall use that  every automorhism of  ${\bf SL}(2,\F_p)$ 
($p$ a prime) is induced by an inner automorphism of the larger group
\begin{equation}
{\bf SL}^{\pm 1}(2,\F_p):=\langle {\bf SL}(2,\F_p), \ 
W:=\begin{pmatrix} 0 & 1 \\ 1 & 0 \end{pmatrix}\rangle.
\end{equation}
See the appendix of \cite{di}  for a proof of this fact. We also use the
following lemma which is easy to prove.

\begin{lem}\label{vd4}

1. Let $H$ be a perfect group. Every automorhism $\psi: H_{[4]}\to H_{[4]}$
satisfies $\psi(H\times H\times \{0\})=H\times H\times \{0\})$.

2. If $H$ is  non-abelian simple finite group, then every automorhism of
   $H\times H$ is of product type.

3. Let $H$ be  ${\bf SL}(2,\F_p)$, where $p$ is prime: then every automorhism of
   $H\times H$ is of product type.

\end{lem}

\Proof
1.: $H\times H\times \{0\}$ is the commutator subgroup.

2.:  the centralizer $ C((x,y))$ of an element $(x,y)$ where 
$ x \neq 1, y \neq 1$ does not map surjectively onto $H$ through either 
of the two product projections, whence every automorphism leaves invariant the
unordered pair of subgroups $\{(H\times \{0\}) ,(\{0\}) \times H) \}.$

3. follows by the same argument used for 2.

\QED

We shall apply the above constructions to obtain some concrete examples.

\begin{prop}\label{vd5}
 Let $p$ be a prime with $p\equiv 3$ mod $4$ and 
 $p\equiv 1$ mod $5$ and consider the group $H:={\bf SL}(2,\F_p)$. 
Then $H_{[4]}$ admits a mixed Beauville structure $u$ such that $\iota(u)$ does
 not lie in the ${\rm A}_\M(H_{[4]})$ orbit of $u$.
\end{prop}
\Proof Set $a_1:=B$, $c_1:=S$ as defined in (\ref{ps1}) and $a_2$, $c_2$ one
of the systems of generators constructed in Proposition \ref{psp3}. That is 
the equations
$$\gamma a_2\gamma^{-1}= a_2^{-1},\quad \gamma c_2\gamma^{-1}= c_2^{-1}$$
are solvable with $\gamma\in{\bf SL}(2,\F_p)$ but not with 
$\gamma\in{\bf SL}(2,\F_p)W$. Set $a:=(a_1,a_2,2)$, $c:=(c_1,c_2,2)$.
By Lemma \ref{vz3} the triple $u:=(H_{[2]},a,c)$ is a 
mixed Beauville structure
on  $H_{[4]}$. The type of $(a,c)$ is $(20,30,5p)$, hence it is strict. 

Suppose that $\iota(u)$ is 
in the ${\rm A}_\M(H_{[4]})$ orbit of $u$. By Corollary \ref{vd1} we have an
automorphism $\psi: H_{[4]}\to H_{[4]}$ with $\psi( H_{[2]})= H_{[2]}$ with
$\psi(a)=a^{-1}$ and $\psi(c)=c^{-1}$. From Lemma \ref{vd4} we get two
elements $\gamma_1,\gamma_2\in {\bf SL}^{\pm 1}(2,\F_p)$ with 
$$\gamma_1 a_1\gamma_1^{-1}= a_1^{-1},\ \gamma_1 c_1\gamma_1^{-1}= c_1^{-1},\
\gamma_2 a_2\gamma_2^{-1}= a_2^{-1},\ \gamma_2 c_2\gamma_2^{-1}= c_2^{-1}.
$$
Since they come from the automorphism $\psi: H_{[4]}\to H_{[4]}$ they have to
lie in the same coset of ${\bf SL}(2,\F_p)$ in 
${\bf SL}^{\pm 1}(2,\F_p)$. This is impossible since 
$\gamma_1 a_1\gamma_1^{-1}= a_1^{-1},\ \gamma_1 c_1\gamma^{-1}= c_1^{-1}$ is
only solvable in the coset  ${\bf SL}(2,\F_p)W$ as a computation shows.
\QED

{\it Proof of Theorem \ref{intro2}}:
We take the mixed Beauville structure from Proposition \ref{vd5} and let $S$
be the corresponding mixed Beauville surface as constructed in Section 4.1.  
\QED

\section{Generating groups by two elements}

\subsection{Symmetric groups  }

In this section we provide a series of intermediate results which
lead to the proof of the above Theorem \ref{sym}. In fact we prove:
\begin{prop}\label{sy11}
Let $n\in \N$ satisfy $n\ge 8$ and $n\equiv 2$ mod $3$, then
$\SSS_n$ has systems of generators $(a,c)$, $(a',c')$ with

1. $\Sigma(a,c) \cap \Sigma(a',c') = \{ 1 \},$

2. there is no $\gamma\in \SSS_n$ with $\gamma a\gamma^{-1}=a^{-1}$ and
  $\gamma c\gamma^{-1}=c^{-1}$.
\end{prop}

\begin{lem}
Let $G$ be the symmetric group $\SSS_n$ in $n \geq 7$ letters,
let  $ a : = (5,4,1) (2,6)$, $ c: = (1,2,3)(4,5,6 .\dots, n)$. There is no
automorphism of $G$ carrying $ a \ra a ^{-1}$, $c \ra c^{-1}$.
\end{lem}

\Proof Since $n \neq 6$, every automorphism of $G$ is an inner one.
If there is a permutation $g$ conjugating $ a $ to $ a ^{-1}$, $c $ to
$ c^{-1}$, $g$ would leave each of the sets $\{1,2,3\}$,
   $\{4,5,.\dots ,n\}$,
$\{1,4,5 \}$, $\{2,6 \}$ invariant. By looking 
at their
intersections we conclude that $g$ leaves the elements 
$1,2,3,6$ fixed
and that the set $\{4,5\}$ is invariant.

Since $g$ 
leaves $1$ fixed and conjugates $a$ to $a^{-1}$, we see moreover
that 
$g$ transposes $4$ and $5$.

But then $g$ conjugates $c$ to 
$(1,2,3)(5,4,6, \dots, n)$ which
is a different permutation than $ 
c^{-1}$.
\QED

\begin{lem}
The two elements $ a : = (5,4,1) (2,6)$, 
$ c: = (1,2,3)(4,5,6, \dots, n)$
generate the symmetric group 
$\SSS_n$ if $n \geq 7$ and $ n \neq 0$ mod $3$.
\end{lem}

\Proof
Let 
$G$ be the subgroup generated by $a,c$.  Then $G$ is generated 
also
by $s, \alpha, T, \gamma$ where $ s := (2,6)$, $ \alpha := 
(5,4,1)$,
$ T := (1,2,3)$, $ \gamma := (4,5,6, \dots n)$,since these 
elements
are powers of $a,c$ and $3$ and $n-3$ are relatively 
prime.

Since $G$ contains a transposition,
it suffices to show that 
it is doubly transitive.

The transitivity of $G$ being obvious, 
since the supports of the
cyclic permutations  $s, \alpha, T, \gamma$ 
have the whole set
$\{1,2, \dots n\}$ as union, let us consider the 
subgroup $ H \subset
G$ which stabilizes $\{ 3\}$. $H$ contains $s, 
\alpha,  \gamma$ and
again these are cyclic permutations such that 
their supports have
as union the set $\{1,2,4,5 \dots n\}$. Thus $G$ 
is doubly transitive,
whence $ G = 
\SSS_n$.
\QED

\begin{oss}

\begin{itemize}
\item
Since $3$ does not 
divide $n$, it follows that
${\rm ord}(c) = 3 (n-3)$,
while ${\rm 
ord}(a)=6$.

\item
We calculate now ${\rm ord}(b)$, recalling 
that
$abc = 1$, whence $b$ is the inverse of $ca$.
Since $ca = 
(1,6,3) (4,2,7,\dots n)$ we have
${\rm ord} (b) = {\rm lcm} (3, 
n-4)$.
\item
Recalling that  $ a' : =  \sigma^{-1}$, $c' : = \tau 
\sigma^2$, where 
$\tau := (1,2)$
and $\sigma: = (1,2, \dots , 
n)$,
it follows immediately that $a', c'$ generate the whole 
symmetric group.
\item
We have ${\rm ord}(b') = {\rm ord}(c'a') = 
{\rm ord} 
(\tau \sigma) = {\rm ord}((2,3,\dots n))
= n-1 $, ${\rm 
ord}(a')  = {\rm ord}( \sigma) = n$.
\item
If $n$ is even, $n=2m$, 
then
$c' = (1,2) (1,3,5,\dots 2m-1)(2,4,6,\dots 2m) $ is the cyclical 
permutation
$ c' = (2,4,\dots, 2m, 1,3,\dots , 2m-1)$ and ${\rm 
ord}(c')=n$.
\item
If $n$ is odd, $n=2m+1$, then
$c' = (1,2) 
(1,3,5,\dots 2m+1,2,4,6,\dots 2m)
  (1,3,5,\dots 2m+1)(2,4,6,\dots 
2m) $ and thus ${\rm ord}(c') = 
m(m+1)$.
\end{itemize}
\end{oss}

\begin{prop}
Let $a,b,c$, 
$a',b',c'$ be as above then
$$   \Sigma(a,c) \cap \Sigma(a',c') = \{ 
1 \} .  $$
\end{prop}

\Proof
We say that a permutation has type 
$(d_1 \leq \dots  \leq d_k)$,
with $d_i \geq 2  \ \forall \ i$, if 
its
cycle decomposition consists of $k$ cycles of respective 
lengths
$d_1, \dots d_k$. We say that the type is monochromatic if 
all
the $d_i$'s are equal, and dichromatic
if the number of distinct 
$d_i$'s is exactly two.
Two permutations are conjugate  to each other 
iff their types
are the same.
We say that a type $(p_1 \leq \dots 
\leq p_r)$ is derived from
$(d_1 \leq \dots  \leq d_k)$ if it is the 
type of a power of
a permutation of type $(d_1 \leq \dots  \leq 
d_k)$.

Therefore we observe that the types of $\Sigma(a,c)$ are 
those derived
from $(2,3)$, $(3,n-3)$, $(3, n-4)$, while those of 
$\Sigma(a',c')$
are those derived from $(n), (n-1)$ for $n$ even, and 
also from
$(m, m+1)$  in case where $ n = 2m+1$ is odd.

We use then 
the following lemma whose proof is straightforward

\begin{lem}
Let 
$g$ be a permutation of type $ (d_1, d_2)$. \\
Then the type of
$ 
g^h$ is the reshuffle of $h_1$-times $d_1/h_1$ and
$h_2$-times 
$d_2/h_2$, where $h_i := {\rm gcd} (d_i, h)$.
Here, reshuffling means 
throwing away all the numbers equal to $1$ and
arranging the others 
in increasing order.

In particular, if the type of $ g^h$ is 
dichromatic, $ (d_1, d_2)$ are
automatically determined. If moreover 
$d_1$, $d_2$ are relatively
prime and the type of $ g^h$ is 
monochromatic,
  then it is derived from type
$d_1$ or from type 
$d_2$.

\end{lem}

For the types in $\Sigma(a',c')$, we get types 
derived from $(n)$, $(n-1)$,
or $(\frac{1}{2} (n-1), \frac{1}{2} 
(n+1))$. The latter come from 
relatively prime numbers,
whence they 
can never equal a type in $\Sigma(a,c)$, derived from the 
pairs
$(2,3)$, $(3,n-4)$  and $(3,n-3$.
The monochromatic types in 
$\Sigma(a,c)$ can only be derived by $(3)$, $(2)$,
$(n-4)$, $(n-3)$, 
since we are assuming that $3$ does neither divide
  $n$ nor 
$n-1$.

\QED

\subsection{Alternating groups}

In this section we 
shall construct certain systems of generators of the
alternating 
groups $\AAA_n$ ($n\in \N$). Our principal tool is the theorem of 
Jordan, see \cite{wi}. This result says that $\langle a,c\rangle= \AAA_n$ 
for any 
pair $a,c\in \AAA_n$ which satifies
\begin{itemize}
\item 
the group  $H:=\langle a,c\rangle$ acts primitively on 
$\{1,\ldots,n\}$,
\item the group $H$ contains a $q$-cycle for a 
prime $q\le n-3.$
\end{itemize}
A further result that we shall need is:

\begin{lem} \label{all1}
For $n\in\N$ with $n\ge 12$ let $U\le \AAA_n$ be a 
doubly transitive group. If $U$ contains a double-transposition then
$U=\AAA_n$. 
\end{lem} 
\Proof
The degree $m(\sigma)$ of a permutation 
$\sigma\in\AAA_n$ is the number of elements moved
by $\sigma$. Let $\sigma\in U$ be a double-transposition. We have
$m(\sigma)=4$. Let $m$ now be the minimal degree taken over all non-trivial
elements of $U$. Since $U$ is also primitive we may apply 
a result of de S\'eguier (see \cite{wi}, page 43) which says that if $m>3$
(in our situation we would have m=$4$) then 
\begin{equation}\label{al1}
n<\frac{m^2}{4} {\rm log}\frac{m}{2}+m\left({\rm log}
\frac{m}{2}+\frac{3}{2}\right).
\end{equation}
For $m=4$ the right-hand side of \ref{al1} is roughly $11.5$. Our
assumptions imply that $m=3$. We then apply Jordan's theorem to reach the
desired conclusion.

\QED

We shall start now to construct the systems of generators required for the
constructions of Beauville surfaces.
We treat permutations as maps which 
act from the left on the set of reference. We also use the notation 
$g^\gamma:=\gamma g\gamma^{-1}$ for the conjugate of an element $g$.
\begin{prop}\label{alp1}
Let $n\in\N$ be even with $n\ge 16$ and let $3\le p\le q\le n-3$ be primes with
$n-q\not\equiv 0$ mod $p$. Then there is a system $(a,c)$ of generators for 
$\AAA_n$ of type $(q,p(n-q),n-p+2)$ such that there is no $\gamma\in \SSS_n$
with $\gamma a\gamma^{-1}= a^{-1}$
and $\gamma c\gamma^{-1}= c^{-1}$.
\end{prop}
\Proof
Set $k:=n-q$ and define
$$a:=(1,2,\ldots,q),\ c:=(q+1,q+2,\ldots,q+k-1,1)
(q+k,p,p-1,\ldots,2).$$
We compute
$$ca=(1,q+k,p,p+1,\ldots,q+k-1)$$
and the statement about the type is clear. We show that there
is no $\gamma\in \SSS_n$ with the above properties. Otherwise,
$\gamma$ would leave invariant the three sets corresponding to the non trivial 
orbits of $a$, respectively $c$, and in particular we would have
$\gamma(1) = 1$, $\gamma (\{2, \dots , p\}) = \{2, \dots , p\}$. 
But then   $\gamma (2) \neq q$, a contradiction.
We set 
$U:=\langle a,c\rangle$ and show that $U=\AAA_n$.

Obviously $U$ is transitive. The stabiliser $V$ of $q+k$ in $U$ contains 
the elements $a,c^p$. It is clear that the subgroup generated by these two
elements is transitive on $\{1,\ldots,n\}\setminus \{ q+k\}$, hence $U$ is
doubly transitive. The group $U$ contains the $q$-cycle $a$, whence we
infer by  Jordan's theorem that $U=\AAA_n$.\QED

For the applications in the previous sections we need:
\begin{prop}\label{alp2}
1. Let $n\in \N$ satisfy $n\ge 16$ with $n\equiv 0 \ {\rm mod}\ 4$ and 
$n\equiv 1 \ {\rm mod}\ 3$. There is a pair $(a,c)$ of generators of 
$\AAA_n$
of type $(2,3,84)$ and an element $\gamma\in \SSS_n\setminus \AAA_n$ with $\gamma 
a\gamma^{-1}= a^{-1}$
and $\gamma c\gamma^{-1}= c^{-1}$.

2. Let $p$ 
be a prime with $p> 5$. Set 
$n:=3p+1$. There is a pair $(a,c)$ of generators of $\AAA_n$
of type $(p,5p,2p+3)$ and $\gamma\in \SSS_n$ with 
$\gamma a\gamma^{-1}= a^{-1}$ and $\gamma c\gamma^{-1}= c^{-1}$. If 
$p\equiv 1$ mod $4$ then $\gamma$ can be chosen in $\SSS_n\setminus \AAA_n$,
if $p\equiv 3$ mod $4$ then $\gamma$ can be chosen in $\AAA_n$.
\end{prop}

\Proof 
1. We take 
as set of reference $\{0,\ldots,n-1\}$ instead of
$\{1,\ldots,n\}$ 
and set 
$$\gamma:=(0)(1)\prod_{i=1}^{\frac{n-4}{6}}(6i-2,6i+1)\cdot 
(2,3)\cdot
\prod_{i=1}^{\frac{n-4}{6}}(6i-1,6i+3)(6i,6i+2),
$$

$$a:=(0,1)\cdot\prod_{i=1}^{\frac{n-4}{6}}(6i-2,6i+1)\cdot 
\prod_{i=1}^{\frac{n-4}{6}}(6i-4,6i-1)\cdot
\prod_{i=1}^{\frac{n-4}{6}}(6i-4,6i-1)^\gamma\cdot
(n-2,n-4),
$$

$$c:=(0)\prod_{i=1}^{\frac{n-1}{3}}(3i-2,3i-1,3i).
$$

Observe that $$\prod_{i=1}^{\frac{n-4}{6}}(6i-4,6i-1)^\gamma =
\prod_{i=2}^{\frac{n-4}{6}}(6i-6,6i+3) \cdot (3,9).$$

We 
have now
$$ca=(0,2,6,13,11,9,1)\cdot(3,7,5)\cdot
\prod_{i=1}^{\frac{n-16}{6}}(6i-2,6i+2,6i+6,6i+13,6i+11,6i+9)
$$
$$
\cdot(n-1,n-12,n-8,n-4)\cdot(n-2,n-6)
$$
We have ${\rm ord}(a)=2$, ${\rm ord}(c)=3$,  ${\rm ord}(ca)=84$ and 
$\gamma a\gamma^{-1}= a^{-1}$
and $\gamma c\gamma^{-1}= c^{-1}$.
To apply the theorem of Jordan we note that $(ca)^{12}$ is a $7$-cycle. 
It remains to show that $H:=\langle a,c\rangle$ acts primitively. 
In figure 1 we exhibit the orbits of $a$ and $c$. In the right hand 
picture we connect two elements of $\{ 0,\ldots,n-1\}$ if they are 
in the 
orbit of $a$, in the left hand picture similarly for 
$c$.

\begin{figure}
\entrymodifiers={++[o][F]}
\SelectTips{cm}{}
{\small
$$
\xymatrix@=10pt{{0} 
& *{} & *{} \\
  {1} \ar@{=}[r] & {2} \ar@{=}[r] & {3} \\
  {4} 
\ar@{=}[r] & {5} \ar@{=}[r] & {6} \\
  {7} \ar@{=}[r] & {8} 
\ar@{=}[r] & {9} \\
  *++{\vdots} & *++{\vdots} & *++{\vdots} \\ 
*++{\vdots} & *++{\vdots} & *++{\vdots} \\
  *++[F-:<10pt>]{n-6} 
\ar@{=}[r] & *++[F-:<10pt>]{n-5}
  \ar@{=}[r] & *++[F-:<10pt>]{n-4} 
\\
  *++[F-:<10pt>]{n-3} \ar@{=}[r]& *++[F-:<10pt>]{n-2} \ar@{=}[r] 
&
  *++[F-:<10pt>]{n-1} } \qquad \qquad
\xymatrix@=10pt{0 \ar@{=}[d] 
& *{} & *{} \\
  {1} & {2} \ar@{=}[d] & {3} \ar @(r,r) @{=} [dd] \\ 
{4} \ar@{=}[d] & {5} & {6} \ar @(r,r) @{=} [dd] \\
  {7} & {8} 
\ar@{=}[d] & {9} \\
  *++{\vdots} & *++{\vdots} & *++{\vdots} \\
*++{\vdots} & *++{\vdots} & *++{\vdots} \\
  *++[F-:<10pt>]{n-6} 
\ar@{=}[d] & *++[F-:<10pt>]{n-5}
  \ar@{=}[u] & *++[F-:<10pt>]{n-4} 
\ar@{=}[dl] \\
  *++[F-:<10pt>]{n-3} & *++[F-:<10pt>]{n-2} & 
*++[F-:<10pt>]{n-1}
  \ar@(r,r)@{=}[uu]}
$$
}
\caption{The orbits of 
$c$, $a$}\label{bild1}
\end{figure}
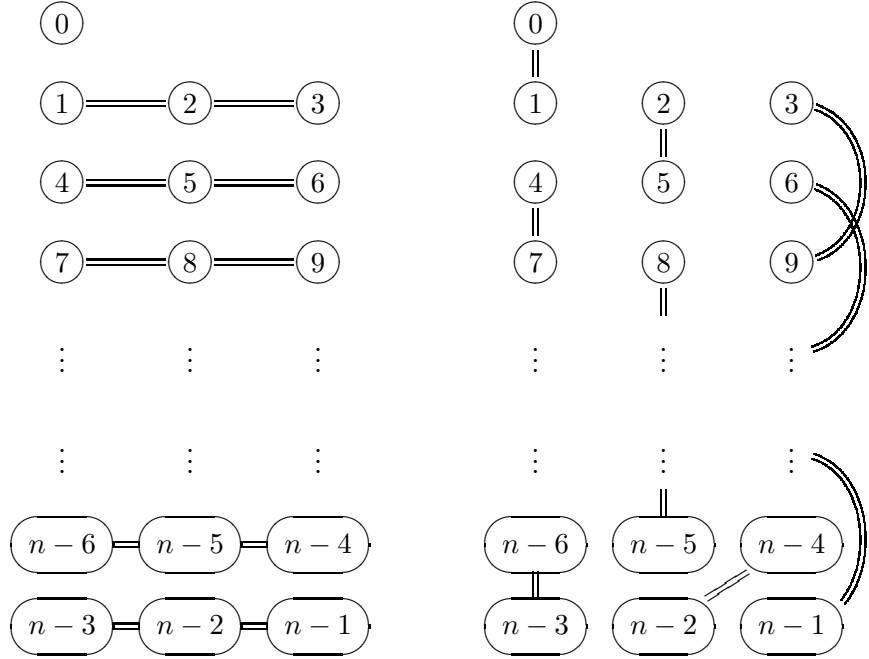

Figure 1 makes it obvious that 
$H$ acts transitively. Since $(ca)^{12}$ is a 
$7$-cycle, the 
second condition of Jordan's theorem is fullfilled. Let 
$H_0$ be the 
stabiliser of $0$. Then $c$ and $(ca)^{7}$ are contained in 
$H_0$. 
In Figure $2$ we show the orbits of the two elements  
$c$ and 
$(ca)^{7}$. The notation is the same as in Figure 1. A glance at 
Figure 2 shows that $H_0$ is transitive on $\{ 1,\ldots,n-1\}$. We 
infer that $H$ is doubly transitive. Again by Theorem 9.6 of 
\cite{wi} the 
group $H$ is 
primitive.

\begin{figure}
\entrymodifiers={++[o][F]}
\SelectTips{cm}{}
{\small
$$
\xymatrix@=10pt{ 
{1} \ar@{=}[r] & {2} \ar@{=}[r] & {3} \\
  {4} \ar@{=}[r] & {5} 
\ar@{=}[r] & {6} \\
  {7} \ar@{=}[r] & {8} \ar@{=}[r] & {9} \\
  {10} 
\ar@{=}[r] & {11} \ar@{=}[r] & {12} \\
  {13} \ar@{=}[r] & {14} 
\ar@{=}[r] & {15} \\
  *++{\vdots} & *++{\vdots} & *++{\vdots} \\
*++{\vdots} & *++{\vdots} & *++{\vdots} \\
  *++[F-:<10pt>]{n-9} 
\ar@{=}[r] & *++[F-:<10pt>]{n-8}
  \ar@{=}[r] & *++[F-:<10pt>]{n-7} 
\\
  *++[F-:<10pt>]{n-6} \ar@{=}[r] & *++[F-:<10pt>]{n-5}
\ar@{=}[r] & *++[F-:<10pt>]{n-4} \\
  *++[F-:<10pt>]{n-3} \ar@{=}[r] 
& *++[F-:<10pt>]{n-2} \ar@{=}[r] &
  *++[F-:<10pt>]{n-1} } \qquad 
\qquad 
\xymatrix@=10pt{
  {1} & {2} & {3} \ar@{=} [dl] \\
  {4} 
\ar@{=}[dr] & {5} \ar@{=}[dl]& {6} \\
  {7} & {8} \ar@{=}[dr] & {9} 
\\
  {10} \ar@{=}[dr] & {11} & {12} \ar@{=}[d] \\
  {13} & {14} 
\ar@{=}[dr] & {15} \ar@{=}[dl] \\
  *++{\vdots} & *++{\vdots} & 
*++{\vdots} \\
  *++{\vdots} & *++{\vdots} & *++{\vdots} \\
*++[F-:<10pt>]{n-9} & *++[F-:<10pt>]{n-8} \ar@{=}[ul] &
*++[F-:<10pt>]{n-7} \ar@{=}[ul] \\
  *++[F-:<10pt>]{n-6} & 
*++[F-:<10pt>]{n-5} \ar@{=}[ur] &
  *++[F-:<10pt>]{n-4} \ar@{=}[ul] 
\\
  *++[F-:<10pt>]{n-3} \ar@{=}[ur] & *++[F-:<10pt>]{n-2} 
\ar@{=}[ul] &
  *++[F-:<10pt>]{n-1} \ar@{=}[u]}
$$
}
\caption{The 
orbits of $c$ and $(ca)^{7}$}\label{bild2}
\end{figure}
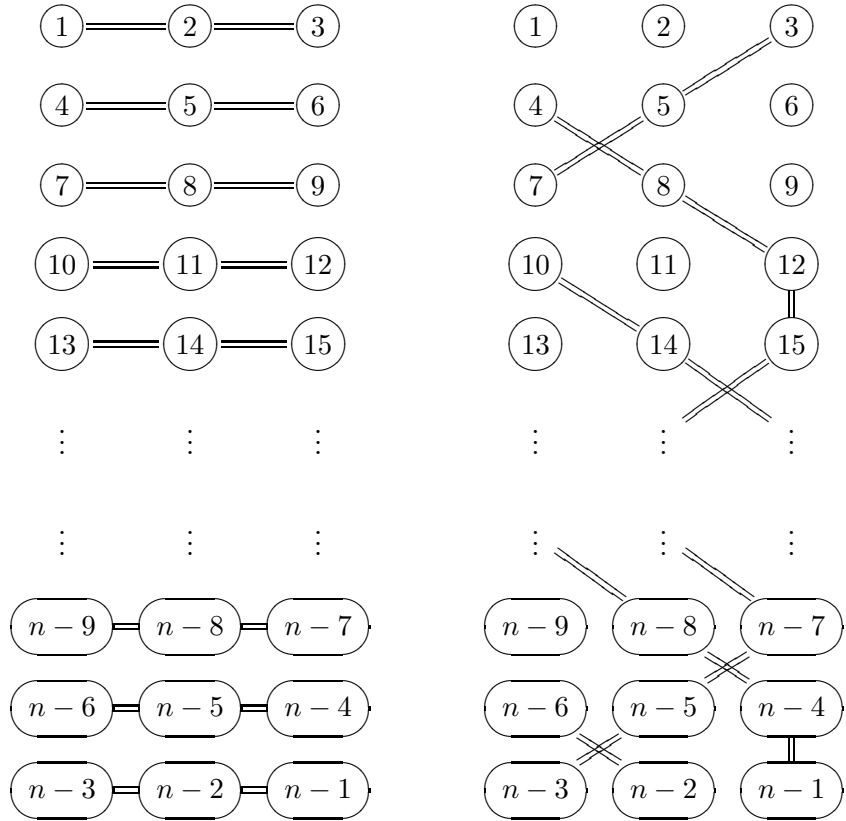

2. Again 
we take as set of reference $\{0,\ldots,n-1\}$ instead 
of
$\{1,\ldots,n\}$ and set
$$\gamma:=(0)(1)(p+1,2p+1)\cdot 
\prod_{i=1}^{\frac{p-1}{2}}(1+i,p+1-i)\cdot
\prod_{i=1}^{p-1}(p+1+i,3p+1-i),
$$
$$a:=(0)(1,2,\ldots,p)(p+1,p+2,\ldots,2p)(2p+1,2p+2,\ldots,3p),
$$
$$c:=(0,p,p-1,p-2,\ldots,3,2)(1,p+1,3p,p+2,2p+1).
$$
We 
have
$$ca:=(2)(3)\ldots(p-1)(0,p,p+1,2p+1,2p+2,\ldots,3p-1,p+2,p+3,\ldots2p,3p,1).
$$

From 
this definition we see that ${\rm ord}(a)=p$, ${\rm ord}(c)=5p$ and 
${\rm ord}(ca)=2p+3$. The formulae   
$\gamma a\gamma^{-1}= a^{-1}$ 
and $\gamma c\gamma^{-1}= c^{-1}$ are also clear.

We verify the 
conditions of Jordan's theorem. First of all $c^5$ is a 
$p$-cycle. 
It remains to show that $H:=\langle a,c\rangle$ acts primitively.
In 
Figure 3, above,  we exhibit the orbits of $a$, respectively the orbits of $c$ below. 
From this it is obvious that $H$ acts 
transitively.

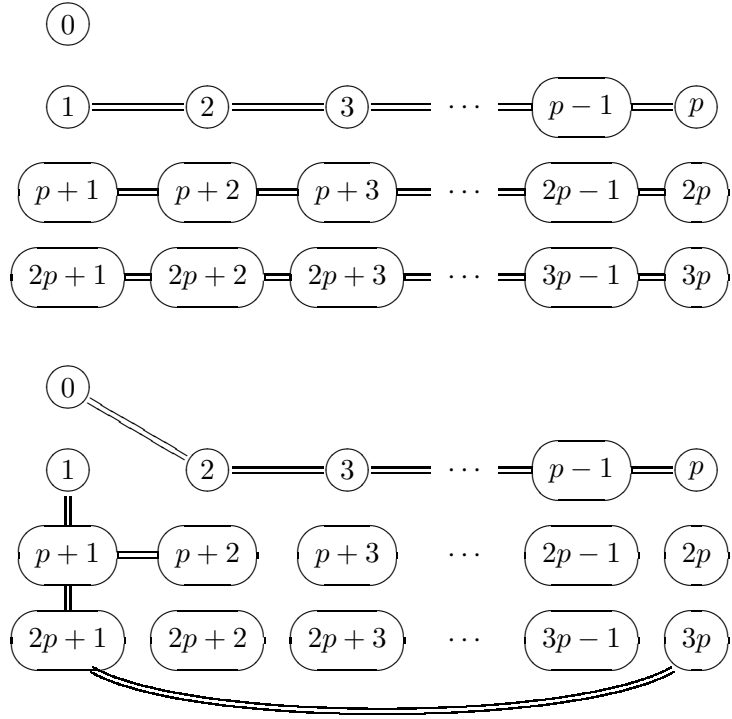
\begin{figure}
\entrymodifiers={++[o][F]}
\SelectTips{cm}{}
{\small
$$
\xymatrix@=10pt{ 
{0} & *{} & *{} & *{} & *{} & *{} \\
  {1} \ar@{=}[r] & {2} 
\ar@{=}[r] & {3} \ar@{=}[r] &
  *++{\cdots} & *++[F-:<10pt>]{p-1} 
\ar@{=}[l] & {p} \ar@{=}[l]\\
  *++[F-:<10pt>]{p+1} \ar@{=}[r] & 
*++[F-:<10pt>]{p+2} \ar@{=}[r] &
  *++[F-:<10pt>]{p+3} \ar@{=}[r] & 
*++{\cdots}
  & *++[F-:<10pt>]{2p-1} \ar@{=}[l] & *++[F-:<10pt>]{2p} 
\ar@{=}[l] \\
  *++[F-:<10pt>]{2p+1} \ar@{=}[r] & 
*++[F-:<10pt>]{2p+2} \ar@{=}[r] &
  *++[F-:<10pt>]{2p+3} \ar@{=}[r] & 
*++{\cdots} & *++[F-:<10pt>]{3p-1}
  \ar@{=}[l] & *++[F-:<10pt>]{3p} 
\ar@{=}[l]
}
$$
\vspace{.5cm}
$$
\xymatrix@=10pt{ {0} \ar@{=}[dr] & 
*{} & *{} & *{} & *{} & *{} \\
  {1} \ar@{=}[d] & {2} \ar@{=}[r] & 
{3} \ar@{=}[r] &
  *++{\cdots} & *++[F-:<10pt>]{p-1} \ar@{=}[l] & {p} 
\ar@{=}[l] \\
  *++[F-:<10pt>]{p+1} \ar@{=}[r] \ar@{=}[d] & 
*++[F-:<10pt>]{p+2} &
  *++[F-:<10pt>]{p+3} & *++{\cdots} & 
*++[F-:<10pt>]{2p-1} &
  *++[F-:<10pt>]{2p} \\
  *++[F-:<10pt>]{2p+1} 
\ar @{=} @(d,d) [rrrrr] & *++[F-:<10pt>]{2p+2} & 
*++[F-:<10pt>]{2p+3} & *++{\cdots} & *++[F-:<10pt>]{3p-1} &
*++[F-:<10pt>]{3p} }
$$
}
\vspace{.5cm}
\caption{The orbits of $a$ 
und $c$}\label{bild3}
\end{figure}

Let now 
$\Delta_1\cup\Delta_2\cup\ldots \cup \Delta_k=\{0,1,\ldots,3p\}$ 
be a block-decomposition for $H$ with $0\in \Delta_1$. 
Note that the 
natural number $k$ satisfies 
$|\Delta_i|k=3p+1$ for all $i$. 
Since 
$a$ is in the stabiliser of $0$ we 
have $|\Delta_1|\in 
\{1,p+1,2p+1,3p+1\}$. 
Since $|\Delta_1|$ divides $3p+1$, we 
infer that $H$ acts primitively.
\QED

Further systems of generators are needed:
\begin{prop}\label{alp3}
Let $n=2k$ be an even natural number with $n\ge 16$. Then there is a system of
generators $(a,c)$ of $\AAA_n$ of type $(2k-3,2k-2,2k-2)$ 
and  $\gamma\in \SSS_n$ such that
$$\gamma a \gamma^{-1}=a^{-1},\quad \gamma c \gamma^{-1}=ac.$$
If $k$ is even then $\gamma$ can be chosen in $\SSS_n\setminus\AAA_n$,
if $k$ is odd then $\gamma$ can be chosen in $\AAA_n$. The system of 
generators $(a,c)$ has the further property that there is no $\delta\in
\SSS_n$ with
$$\delta a \delta^{-1}=a^{-1}\ {\it and}\ \delta c \delta^{-1}=c^{-1}\ \
{\it or}\ \ 
\delta a \delta^{-1}=c^{-1}\ {\it and}\ \delta c \delta^{-1}=a^{-1}.$$
\end{prop}
\Proof We set 
$$a:=(1,2,\ldots,2k-4,2k-3),$$ 
$$d:=(1,2,3)(2k-3,2k-4,2k-5)(k-1,2k-1)(2k-2,k-2,2k,k)
$$
and
$$\alpha:=(1,2k-3)(2,2k-4)(3,2k-5)\ldots (k-2,k)\cdot (2k-2,2k).$$
The following are clear:

\begin{itemize}
\item $d^6$ is a double-transposition,
\item $\alpha$ is in $\AAA_n$ if $k$ is odd, and in $\SSS_n\setminus\AAA_n$
if $k$ is even, 
\item $\alpha a \alpha^{-1}=a^{-1}$ and $\alpha d \alpha^{-1}=d$,
\item there is no $\delta\in \SSS_n$ with
$\delta a \delta^{-1}=a$ and $\delta d \delta^{-1}=d^{-1}$.
\end{itemize}
For the last item note that $a$ has $\{ 2k-2,2k-1,2k\}$ as its set of fixed
points. A $\delta$ with the above property would have to stabilise this
set. The element $d$ interchanges $k-1$ and $2k-1$ hence both these elements
have to be fixed by $\delta$. The condition $\delta a \delta^{-1}=a$ implies
then that $\delta$ acts as the identity on $\{1,2, \dots, 2k-3\}$, in particular as
the identity on the subset  $\{1,2, 3\}$, contradicting $\delta d \delta^{-1} = d^{-1}$.  
Set
$$c:=d
a^{k-2}=(1,2k-1,k-1,2k-4,k-3,2k-3,2k,k,$$
$$
2,2k-2,k-2,2k-5,k-4,2k-6,k-5,2k-7,\ldots,5,k+3,4,k+2)(3,k+1).$$
Note that ord$(c)=2k-2$. We find
$$\alpha c \alpha^{-1}=\alpha da^{k-2} \alpha^{-1}=da^{-k+2}=
ca^{-2k+4}=ca,$$
whence also $ord (ca) = 2k-2$.
Set now $\gamma:=a \alpha$, so that $\gamma$ clearly satisfies the required properties.

Set $U:=\langle a,c\rangle.$ We shall now show that $U=\AAA_n$. Clearly, $U$
is a transitive group. Let $V\le U$ the stabiliser of $2k-1$. The subgroup $V$
contains
$$a,\ d^2,\ da^{k-2}d a^{4-k} d.$$
From the definitions it is clear that these elements generate a group which is
transitive on $\{1,\ldots,2k\}\setminus \{2k-1\}$. Hence $U$ is doubly
transitive and contains a double transposition. From Lemma \ref{all1} we infer
that $U=\AAA_n$. 

The last property follows from the above items since in the first case we would have
$ a \ra a^{-1}, d \ra (Int_a)^{-(k-2)}(d^{-1}) $, and composing with 
$Int_{\alpha \cdot a^{k-2}}$ we contradict the third item.

Whereas, in the second case, just observe that $a$ and $c$ have different order.

\QED

\subsection{SL(2) and PSL(2) over 
finite fields}

In this section we give systems of generators 
consisting of two elements of
the respective groups ${\bf SL}(2,\F_p)$ and ${\bf 
PSL}(2,\F_p)$ which will allow us to
construct certain Beauville 
structures on them. 

If $p$ is a prime we denote by $\F_p$, $\F_{p^2}$ 
 the fields with $p$, respectively $p^2$
elements and by $\F_p^*$, $\F_{p^2}^*$ 
 the corresponding multiplicative
groups. We let 
$$ 
N_{\F_{p^2}/\F_{p}} :\F_{p^2}^*\to \F_{p}^*$$
be the norm map. We 
also
introduce the matrices 
\begin{equation}\label{ps1}
B:=\begin{pmatrix} 0 & 1 \\ -1 & 0 
\end{pmatrix},\
S:=\begin{pmatrix} 0 & -1 \\ 1 & 1 
\end{pmatrix},\
T:=\begin{pmatrix} 1 & 1 \\ 0 & 1 
\end{pmatrix},\
W:=\begin{pmatrix} 0 & 1 \\ 1 & 0 \end{pmatrix}
\end{equation}
in ${\bf GL}(2,\F_p)$. For $\lambda\in \F_p$ with $\lambda \ne 0$ 
and
$k\in\F_p$ we define
$$ D(\lambda):=\begin{pmatrix} \lambda & 0 
\\ 0 & \lambda^{-1} \end{pmatrix},
\quad M(k):=\begin{pmatrix} 0 & 1 
\\ -1 & k \end{pmatrix}.$$ 
We have 
$$B^4=S^6=W^2=1,\ T=BS,\ 
WBW^{-1}=B^{-1},\ WSW^{-1}=S^{-1}.$$
The matrices $(B, S)$ form a 
system of generators of ${\bf SL}(2,\F_p)$ of
type $(4,6,p)$. Their 
images in  ${\bf PSL}(2,\F_p)$ form a system of
generators of type 
$(2,3,p)$.

\begin{prop}\label{psp1}
Let $p$ be an odd prime and let $q\ge 5$ 
be a prime with $q | p-1$ and let $\lambda\in \F_p^*$ be of
order $q$: then $\lambda+\lambda^{-1}-2\ne 0$, 
$\lambda^2+\lambda^{-2}-2\ne 0$ and 
$\lambda+\lambda^{-1}-\lambda^2-\lambda^{-2}\ne 0$. 
Set
\begin{equation} \label{g1} g:=\begin{pmatrix} 1 & b \\ 1 & d 
\end{pmatrix}\
\ {\rm with}\ 
b:=
\frac{\lambda+\lambda^{-1}-\lambda^2-\lambda^{-2}}{\lambda^2+\lambda^{-2}-2},\ 
d:=\frac{\lambda+\lambda^{-1}-2}{\lambda^2+\lambda^{-2}-2}.
\end{equation}
Then 
\begin{equation}\label{ps2} 
D(\lambda),\quad  g D(\lambda) g^{-1}
\end{equation} 
form a system of generators of ${\bf SL}(2,\F_p)$
of type $(q,q,q)$. Set further
\begin{equation} \label{psh3} 
e(\lambda):=\frac{2-\lambda-\lambda^{-1}}{\lambda+\lambda^{-1}-
\lambda^2-\lambda^{-2}}.
\end{equation} 
There exists $\gamma\in {\bf SL}(2,\F_p)$ with 
$\gamma D(\lambda)\gamma^{-1}=D(\lambda)^{-1}$ and 
$\gamma g D(\lambda) g^{-1} \gamma^{-1}=gD(\lambda)^{-1}g^{-1}$
if and only if $e(\lambda)$ is a square in $\F_p$. 
There is  $\gamma\in {\bf SL}(2,\F_p)W$ satisfying these conditions
if and only if $-e(\lambda)$ is a square in $\F_p$.
\end{prop}

\Proof  Let us prove that
$\lambda+\lambda^{-1}-2= 0$ leads to a contradiction. In fact, multiplying 
by $\lambda$ we find 
$(\lambda-1)^2=0$ which is impossible since $\lambda$ has order $q$. The other
two cases are treated similarly.

We see immediately that the determinant of $g$ is equal to $1$, furthermore
an easy computation shows that
the trace of
$h:=D(\lambda) gD(\lambda)g^{-1}$ is 
$\lambda+\lambda^{-1}$, hence $h$ has also order $q$. 
Notice further 
that the subgroup $H$ generated by $D(\lambda)$ and 
$gD(\lambda)g^{-1}$ cannot be solvable because these two elements 
have no
common fixpoint in the action on $\PP^1_{\F_p}$.
From the list of
isomorphism classes of subgroups of 
${\bf PSL}(2,\F_p)$ given in 
\cite{hu} (Hauptsatz 8.27, page 213).  
we find that $H$ has to be ${\bf 
SL}(2,\F_p)$ if $q\ne 5$. For $q=5$ these 
simple arguments with 
subgroup orders leave the possibility that $H$ be isomorphic to the
binary icosahedral group
$2\cdot\AAA_5\le {\bf SL}(2,\F_p) $. But this group 
does not have a system of generators 
of type $(5,5,5)$, as it is easily seen by computer calculation.

The assertion about the simultaneous conjugacy of 
$D(\lambda)$, $D(\lambda)^{-1}$ and 
$g D(\lambda) g^{-1}$, $g D(\lambda)^{-1} g^{-1}$ follows by a straightforward
computation. In fact we just take a matrix $X$ with indeterminate entries
and write down the $9$ equations resulting from      
$ X D(\lambda)=D(\lambda)^{-1} X$ and 
$X g D(\lambda) g^{-1} =gD(\lambda)^{-1}g^{-1} X$, ${\rm Det}(X)=\pm 1$ 
and the statement follows by a small manipulation of them.

\QED

\begin{prop}\label{psp2}
Let $p, q$ be  odd primes such
that $q\ge 5$ and $q | p+1$ . Let $\lambda\in 
\F_{p^2}^*$   be of order 
$q$ (thus $ N_{\F_{p^2}/\F_{p}}(\lambda)=1$). 
Consider its trace $k:=\lambda+\lambda^{-1} \in \F_p$:
then there exists 
$g\in {\bf SL}(2,\F_p)$ such that 
\begin{equation} \label{ps3} 
M(k),\quad  g M(k) g^{-1}
\end{equation} 
form a system of generators 
of type $(q,q,q)$.
\end{prop}

\Proof Choose $k$ as indicated and notice that $k\ne 1$, $k\ne 2$, as we already saw. 

Let us now change our perspective and let $R:=\Z[r,s,t]$ treating $r=k$ as a
variable. Set
$$x:=\begin{pmatrix} 0 & 1 \\ -1 & r 
\end{pmatrix},\ g:=\begin{pmatrix} 1 & s \\ t & 1+st 
\end{pmatrix},\ y:= gxg^{-1},\ z:=x \cdot y\in {\bf SL}(2,R).$$
The reduction into ${\bf SL}(2,\F_p)$ of 
$x$ and of $y$ have order $q$ for any choice of
$s,t$ and $r=k$. We want $z=xy$ also to have order $q$. This happens if 
the trace of the reduction into ${\bf SL}(2,\F_p)$ of $z$ is equal to
$k$. This follows from ${\rm Tr}(z)=r$ which is equivalent to the equation
\begin{equation}\label{psh1}
-s^2t^2 + s^2tr - s^2 - st^2r + str^2 - 2st - t^2 + r^2 - r - 2=0.
\end{equation}
We write $C_{p,k}$ for 
the plane affine curve obtained  from (\ref{psh1}) by setting $r=k$ and
reducing modulo $p$. Furthermore let $\PP C_{p,k}$ be 
the projective closure of $C_{p,k}$ (with respect to $s,\, t$) and 
$\PP C_{p,k}^\infty$ its set of points at infinity.

These are immediately seen to be the two points with $s=0, t=1$, respectively
$t=0, s=1$, and in these two points at $\infty$ the curve has two ordinary double points. 
 
By a computation using a 
Gr\"obner-routine over $\Z$ (possible in MAGMA or SINGULAR)
one can also verify that 
the affine curve $C_{p,k}$ is smooth for every $p\ge 11$.

To check this  fact we projectivise (\ref{psh1}), compute derivatives and
analise the ideal in $R$ generated by these homogeneous polynomials. In this
step we use $k\ne 1$, $k\ne 2$. 

Let now $\tilde C_{p,k}\to \PP C_{p,k}$ be a non-singular model of 
$\PP C_{p,k}$. We conclude from the above analysis of the singularities of
$\PP C_{p,k}$ and from B\'ezout's
theorem that the non-singular curve $\tilde C_{p,k}$ is absolutely irreducible and has
genus $1$. We
may then apply the Hasse-Weil estimate (cf. e.g. the textbook \cite{har}, V 1.10, page
368) to obtain:
$$  ||\tilde C_{p,k}(\F_p)|-p-1|\le 2\sqrt{p}$$
and since there are at most two $\F_p$-points over every singular point of
$\PP C_{p,k}$ we get at worst
$$  ||C_{p,k}(\F_p)|-p+3|\le 2\sqrt{p}.$$
Hence we have $C_{p,k}(\F_p)\ne \emptyset$ for $p\ge 11$.

Let now $(s,t)$ be in $C_{p,k}(\F_p)$ and $x,y,z$ the correspondinding 
matrices defined  above. Notice that all $3$ of them have order $q$.
It can be checked, again by a Gr\"obner-routine over $\Z$ that we have
$y\ne \pm x$ and $y\ne \pm x^{-1}$. Assume that $q>5$. 
A glance at the sugroups of 
${\bf PSL}(2,\F_p)$ (\cite{hu}) shows that the subgroup 
generated by $x,y$ could only be
cyclic which is impossible by the remarks just made. If $q=5$ we conclude by
observing again that $2\cdot\AAA_5\le {\bf SL}(2,\F_p)$ does not have a 
system of generators of type $(5,5,5)$.

\QED

In order to use Proposition \ref{psp1} effectively for our problems we would 
have to show that the invariant $e$ from (\ref{psh3}) takes both square and
nonsquare values as $\lambda$ varies over all elements of order $q$. 
This leads to a difficult problem about exponential sums which we could not
resolve. In case $q=5$ we found the following way to treat the problem by a
simple trick.

\begin{prop} \label{psp3}
Let $p$ be a prime with $p\equiv 3$ mod $4$ and 
 $p\equiv 1$ mod $5$. Then the group ${\bf SL}(2,\F_p)$ has a system of
 generators $(a,c)$ of type $(5,5,5)$ such that the equations 
\begin{equation}\label{pshi3}
\gamma a\gamma^{-1}= a^{-1},\quad \gamma c\gamma^{-1}= c^{-1}
\end{equation}
are solvable with $\gamma\in{\bf SL}(2,\F_p)$ but not with 
$\gamma\in{\bf SL}(2,\F_p)W$. The same group has another system of
 generators $(a,c)$ such that (\ref{pshi3}) is solvable in  
${\bf SL}(2,\F_p)W$ but not in ${\bf SL}(2,\F_p)$.
\end{prop}
\Proof Take $\lambda\in\F_p$ with $\lambda^5=1$, $\lambda\ne 1$ and consider
the system of generators given in (\ref{ps2}). 
Since  $p\equiv 3$ mod $4$ the
number $-1$ is not a square in $\F_p$. Suppose that the invariant $e(\lambda)$
is a square in $\F_p$ whence(\ref{pshi3}) is solvable in $ {\bf SL}(2,\F_p)$ and
not in ${\bf SL}(2,\F_p)W$ (see Proposition \ref{psp1}). Then we are done. Suppose 
instead that the invariant $e(\lambda)$ is not a square in $\F_p$. We replace $\lambda$ by
$\lambda^2$ and find by a small computation that $e(\lambda)=-e(\lambda^2)$ up
to squares. In this place we use $\lambda^5 = 1$,whence $(\lambda^2)^2 = \lambda^{-1}$ 
and the denominator simply changes sign as we replace $\lambda$ by $\lambda^2 $.
Notice also that
$$2-\lambda-\lambda^{-1}=-(\mu-\mu^{-1})^2 \qquad (\mu^2=\lambda)$$
is never  a square. We infer that $e(\lambda^2)$ is a square and proceed as
before. 

The second statement is proved similarly.\QED  

\subsection{Other groups and more generators}

In this subsection we report on computer experiments related to the existence
of unmixed or mixed Beauville structures on finite groups. We also try to
formulate some conjectures concerning these questions.

We have paid special attention to unmixed Beauville structures on finite
nonabelian simple groups. 

The smallest of these groups is $\AAA_5\cong  
{\bf PSL}(2,\F_5)$. This group cannot have an unmixed Beauville structure. On
the one hand it has only elements of orders $1,2,3,5$. It is not solvable
hence it cannot be a quotient group of one of the euclidean triangle groups
(see Section 6). This implies that any normalised system of generators has
type $(n,m,5)$ with $n,m\in \{ 2,3,5\}$. Finally we note that, by Sylow's
theorem, all subgroups of order $5$ are conjugate. 

There are $47$ finite simple nonabelian groups of order  $ \le 50000$. 
By computer calculations we have found  unmixed Beauville structures on all of
them with the exception of $\AAA_5$. This and the results of Section 3.2 
leads us to:

{\bf Conjecture 1:} {\it All finite simple nonabelian groups except $\AAA_5$ 
admit an unmixed Beauville structure.}

We have also checked this conjecture for some bigger simple groups like the  
Mathieu groups ${\bf M}12, \,{\bf M}22 $ and also matrix groups of size bigger then
$2$. Furthermore we have proved:
\begin{prop}
Let $p$ be an odd prime: then the Suzuki group ${\bf Suz}(2^p)$ has an unmixed
Beauville structure.
\end{prop}
In the proof, which is not included here, 
we use in an essential way that the Suzuki groups 
${\bf Suz}(2^p)$ are minimally simple, that is have only solvable proper
subgroups. For the Suzuki groups see \cite{hub}.

Let us call a type $(r,s,t)\in \N$ {\it hyperbolic} if
$$\frac{1}{r}+\frac{1}{s}+\frac{1}{t} < 1.$$
In this case the triangle group $T(r,s,t)$ is hyperbolic.
From our studies also the following looks suggestive: 

{\bf Conjecture 2:} {\it Let  $(r,s,t)$,  $(r',s',t')$ be two hyperbolic
  types. Then almost all alternating groups $\AAA_n$ have an unmixed Beauville
  structure $v=(a_1,c_1;a_2,c_2)$ where $(a_1,c_1)$ has type  $(r,s,t)$ and 
 $(a_2,c_2)$ has type  $(r',s',t')$.}

\medskip
Let us call an  unmixed Beauville structure $v=(a_1,c_1;a_2,c_2)$ on the
finite group $G$ {\it strongly real} if there are $\delta_1,\delta_2\in G$ 
and $\psi\in {\rm Aut}(G)$ with 
\begin{equation}
(\delta_1\psi(a_1)\delta_1^{-1},\delta_1\psi(c_1)\delta_1^{-1};
\delta_2\psi(a_2)\delta_2^{-1},\delta_2\psi(c_2)\delta_2^{-1})=
(a_1^{-1},c_1^{-1};a_2^{-1},c_2^{-1}).
\end{equation}
If the unmixed Beauville structure $v$ is strongly real then the associated
surface $S(v)$ is real.

There are $18$ finite simple nonabelian groups of order  $\le 15000$. 
By computer calculations we have found  
strongly unmixed Beauville structures on all of
them with the exceptions of 
$\AAA_5$, ${\bf PSL}(2,\F_7)$, $\AAA_6$, $\AAA_7$,
${\bf PSL}(3,\F_3)$, ${\bf U}(3,3)$ and the Mathieu group ${\bf M}11$.  
The alternating group $\AAA_8$ however has such a structure.
This and the results of Section 3.3 
leads us to:

{\bf Conjecture 3:} {\it All but finitely many finite simple groups have 
a strongly real unmixed Beauville structure.}

Conjectures 1,2,3 are variations of a conjecture of Higman saying that
every hyperbolic triangle group surjects onto almost all alternating
groups. This conjecture was resolved positively in \cite{ev} where a
related discussion can be found.  

\medskip

We were unable to find finite $2$- or $3$-groups  having an unmixed
Beauville structure. For $p\ge 5$ our construction (\ref{x}) gives plenty of
examples of $p$-groups having an unmixed
Beauville structure.

\medskip

Finally we report now on two general facts that we have found during our
investigations. These are useful in the quest of finding Beauville structures
on finte groups.

Using the methods used the proofs of Propositions \ref{psp1}, \ref{psp2} the
following can be proved:
\begin{prop} Let $p$ be an odd prime.
 
1. Let $q_2>q_1 \ge 5$ be primes with $q_1q_2 | 
p-1$ and 
   let $\lambda_1,\lambda_2\in \F_p^*$ be of
   respective orders 
$q_1$ and $q_2$. Then there is an element $g\in {\bf SL}(2,\F_p)$ such that 
\begin{equation} \label{ps4} 
D(\lambda_1),\quad  g D(\lambda_2) g^{-1}
\end{equation} 
form a system of generators of type 
$(q_1,q_2,q_1q_2)$.

2. Let $q_2>q_1\ge 5$ be primes with $q_1q_2 | 
p+1$ and let 
   $\lambda_{1,2}\in \F_{p^2}^*$ with $ 
N_{\F_{p^2}/\F_{p}}(\lambda_{1,2})=1$ 
   be of respective orders $q_1$ and 
$q_2$. Then their traces
   $k_{1,2}:=\lambda_{1,2}+\lambda^{-1}_{1,2}$ are in 
$\F_p$ and there is an element 
   $g\in {\bf SL}(2,\F_p)$ such that 
\begin{equation} \label{ps5} M(k_1),\quad  g M(k_2) 
g^{-1}
\end{equation}
form a system of generators of type 
$(q_1,q_2,q_1q_2)$.
\end{prop}

Surfaces $S$ which are not real but still are biholomorphic to their conjugate
$\bar S$ are somewhat difficult to find. Our Theorem \ref{intro3} gives
examples using the alternating groups. We also have found:

\begin{prop}\label{r1} Let $p$ be an odd prime and assume that there is a
prime $q\ge 7$ dividing $p+1$ such that $q$ is not a square modulo $p$ : then
there is  an unmixed Beauville surface $S$ with group $G={\bf SL}(2,\F_p)$
which is biholomorphic to the complex conjugate surface $\bar{S}$ but is not
real.
\end{prop}

For the proof we turn the conditions into polynomial equations and 
polynomial inequalities (as in Propositions \ref{psp1}, \ref{psp2}) and 
then use arithmetic algebraic geometry over finite fields (in a more
subtle way) as before. We do not include this here.

\begin{oss} First examples of primes $p$ satisfying the conditions of
  Proposition 
   \ref{r1} are $p=13$ with $q=7$, $p=37$ with $q=19$ and $p=41$ with $q=7$.
Let $p,\, q$ be odd primes. The
law of quadratic reciprocity implies that the conditions of 
Proposition \ref{r1}
are equivalent to $q\equiv 3$ mod $4$, $p\equiv 1$ mod $4$ and $p\equiv -1$
mod $q$. Dirichlet's theorem on primes in arithmetic progressions implies
that there are infinitely many such pairs $(p,\, q)$.
\end{oss}

\section{The wall paper 
groups}\label{wp}

In this section we analyse the finite quotients of 
the triangular groups
$$T(3,3,3),\qquad T(2,4,4), \qquad 
T(2,3,6).$$
and we will show that they do not admit any unmixed Beauville structure.
We shall give two proofs of this fact, a "geometric" one, and the other
in the taste of combinatorial group theory.

These are groups of motions of the euclidean plane, in 
fact in the classical
classification they are the groups {\bf p3}, 
{\bf p4},  
{\bf p6}. Each of them contains a normal subgroup $N$ 
isomorphic to $\Z^2$
with finite quotient. 

In fact, let $T$ be such a triangle group: then $T$ admits a maximal surjective
homomorphism onto a cyclic group $C_d$ of order $d$. Here, $d$ is 
respectively equal
to $3,4,6$, and the three generators map to elements of $C_d$ whose 
order equals their
order in $T$.

It follows that the covering corresponding to $T$ is the universal cover of
the compact Riemann surface $E$ corresponding to the surjection onto $C_d$,
and one sees immediately two things:

1)  $E$ is an elliptic curve because $\mu(a,c) = 1$

2) $E$ has multiplication
by the group $\mu_d \cong \Z / d\Z$ of $d$-roots of unity.

Letting $\omega = exp ( 2/3  \pi i )$, we see that

\begin{itemize}
\item
$T(3,3,3)$ is the group of affine transformations of $\C$ of the form
$$ g(z)  = \omega^j  z + \eta \ , {\rm for} \  j \in \Z / 3\Z , \eta 
\in \Lambda_{\omega}
: =  \Z \oplus \Z \omega $$
\item
$T(2,4,4)$ is the group of affine transformations of $\C$ of the form
$$ g(z)  = i^j  z + \eta \ , {\rm for} \ j \in \Z / 4\Z, 
\eta  \in \Lambda_{i} :
=  \Z \oplus \Z i $$
\item
$T(2,3,6)$ is the group of affine transformations of $\C$ of the form
$$ g(z)  = (- \omega)^j  z + \eta \ , {\rm for} \  j \in \Z / 6\Z , \eta  \in
\Lambda_{\omega} : =  \Z \oplus \Z \omega .$$

\end{itemize}

\begin{oss}
Using the above affine representation, we see that $N$ is the normal 
subgroup of
translations, i.e. of the transformations which have no fixed point on $\C$.

Moreover, if an element $g \in T - N$, then the linear part of $g$ is in
$\mu_d - \{ 1\}$, and $g$ has a unique fixed point $p_g$ in $\C$.
An immediate calculation shows that indeed this fixed point $p_g$ 
lies in the lattice
$\Lambda$, and we obtain in this way that the conjugacy classes of elements
  $g \in T - N$ are exactly given by their linear parts, so they are 
in bijection with the
elements of $\mu_d - \{ 1\}$.
\end{oss}

Let now $G = T / M $ be a non trivial finite quotient
group of $T$:
then $G$ admits a maximal surjective homomorphism onto a cyclic group $C'$
of order $d$, where $ d  \in \{ 2,3,4,6\}$. Assume that there is an element
$ g \in T-N$ which lies in the kernel of the composite homomorphism: 
then the whole
conjugacy class  of $g$ is in the kernel.
Since all transformations in the $N$ - coset of $g$  are in the 
conjugacy class,
it follows that $N$ is in the kernel and $G$ is cyclic, whence 
isomorphic to $C'$.

In the case where $C'$ is isomorphic to $C$, we get that
$G$ is a semidirect product
$ G = K \rtimes C$, where $ K = N / N \cap M$, and the
action of $C$ on $K$ is induced by the one of $C$ on $N$. We have thus shown:

\begin{prop}
Let $G$ be a non trivial finite quotient of a triangle group $ T = 
T(3,3,3),\ {\rm or} \
T(2,4,4), \  {\rm or} \ T(2,3,6)$. Then there is a maximal surjective 
homomorphism
of $G$ onto a cyclic group $C_d$ of order $ d \leq 6$.

If moreover $G$ is not  isomorphic to $C$, then $d=3$ for $T(3,3,3)$,
for $T(2,4,4)$ $d=4$ ,  $d= 6$ for $T(2,3,6)$,
and $G$ is a semidirect product $ G = K \rtimes C$, where the action of
$C$ is induced by the one of $C$ on $N$.
In particular, let $a_1, c_1$ and $a_2, c_2$ by
two systems of generators of $G$: then $|\Sigma(a_1,c_1) \cap 
\Sigma(a_2,c_2)|\neq
\emptyset$.
\end{prop}

\Proof
Just observe that two elements which have the same image in $C - \{ 0 
\}$ belong to the
same conjugacy class by our previous remarks. The rest follows rightaway.
\QED

We give now an alternative proof by purely group theoretical arguments.

In case of $T(3,3,3)$ we have an 
isomorphism of finitely presented groups 
$$\langle a,c \, |\, a^3, 
c^3, (ac)^3\rangle \cong
\langle x,y,r \,|\, 
[x,y],r^3,rxr^{-1}=y,ryr^{-1}=x^{-1}y^{-1}\rangle$$
given by 
$x=ca^{-1}$, $y=cac$, $r=a$. We set $N_3:=\langle x, y\rangle$. 
The 
second presentation shows that $\Gamma(3,3,3)$ is 
isomorphic to the 
split extension of $N_3\cong \Z^2$ by the cyclic group (of
order 3) 
generated by $r$. We have

\begin{prop}
Let $L$ be a normal subgroup 
of finite index in $T(3,3,3)$. If $L\ne
T(3,3,3)$ then $L\le N_3$ and 
$G:= T(3,3,3)/ L$ 
is isomorphic to the split extension of
a finite 
abelian group $N$ by a cyclic group of order $3$. The only possible 

types for a two generator system of $G$ are (up to permutation) 
$(3,3,3)$ and 
$(3,3,l)$ for some divisor $l$ of $ |N|$. Let $a_1, 
c_1$ and $a_2, c_2$ by two
systems of generators of $G$ then 
$|\Sigma(a_1,c_1) \cap \Sigma(a_2,c_2)|\ge
3$.
\end{prop}

\Proof An 
obvious computation shows that the
normal closure of any element 
$g=ur$ ($u\in N_3$) contains $N_3$ and hence is
equal to $T(3,3,3)$. 
This proves the first statement. Let now $L\le N_3$
and let $a_1, 
a_2$ generate $G= T(3,3,3)/ L$ then at least one of the
cosets  $a_1, 
a_2$ must contain an element of the form $g=ur^{\pm 1}$ ($u\in
N_3$). 
By rearrangement both cosets contain an element of this type. 
A 
computation shows that every element has order exactly $3$ 
in
$T(3,3,3)$. This shows he statement about the types. Let 
$g=ur^{\pm 1}$
be as above and let $S$ be the union of the conjugates 
of the cyclic group
generated by $g$ in $\Gamma(3,3,3)$. It is clear 
that $S$ contains either $xr$
or $r$ or both these elements. 
\QED

In case of $T(2,4,4)$ we have an isomorphism of finitely 
presented groups 
$$\langle a,c \, |\, a^2, c^4, (ac)^4\rangle 
\cong
\langle x,y,r \,|\, 
[x,y],r^4,rxr^{-1}=y,ryr^{-1}=y^{-1}\rangle$$
given by $x=ac^2$, 
$y=cac$, $r=c$. We set $N_4:=\langle x, y\rangle$. 
The second 
presentation shows that $\Gamma(2,4,4)$ is 
isomorphic to the split 
extension of $N_4\cong \Z^2$ by the cyclic group (of
order 4) 
generated by $r$. We have

\begin{prop}
Let $L$ be a normal subgroup 
of finite index in $T(2,4,4)$. If the index
of $L$ in $T(2,4,4)$ is 
$\ge 16$ then  $L\le N_4$ and 
$G:= T(2,4,4)/ L$ is isomorphic to the 
split extension of
a finite abelian group $N$ by a cyclic group of 
order $4$. The only possible 
types for a two generator system of $G$ 
are (up to permutation) $(2,4,4)$ and 
$(4,4,l)$ for some divisor $l$ 
of $ |N|$. Let $a_1, c_1$ and $a_2, c_2$ by two
systems of generators 
of $G$ then $|\Sigma(a_1,c_1) \cap 
\Sigma(a_2,c_2)|\ge
2$.
\end{prop}
The proof is analogous to the 
first Proposition of this section.

In case of $T(2,3,6)$ we have an 
isomorphism of finitely presented groups 
$$\langle a,c \, |\, a^2, 
c^3, (ac)^6\rangle \cong
\langle x,y,r \,|\, 
[x,y],r^6,rxr^{-1}=y^{-1}x,ryr^{-1}=x\rangle$$
given by 
$x=cac^{-1}a$, $y=c^{-1}aca$, $r=ac$. 
We set $N_6:=\langle x, 
y\rangle$. 
The second presentation shows that $\Gamma(2,3,6)$ is 
isomorphic to the split extension of $N_6\cong \Z^2$ by the cyclic 
group (of
order 6) generated by $r$. We have

\begin{prop}
Let $L$ be 
a normal subgroup of finite index in $\Gamma(2,3,6)$. If the index
of 
$L$ in $T(2,3,6)$ is $\ge 24$ then  $L\le N_6$ and 
$G:= 
\Gamma(2,3,6)/ L$ is isomorphic to the split extension of
a finite 
abelian group $N$ by a cyclic group of order $6$. The only possible 
types for a two generator system of $G$ are (up to permutation) 
$(2,3,6)$ and 
$(6,6,l)$ for some divisor $l$ of $ |N|$. Let $a_1, 
c_1$ and $a_2, c_2$ by two
systems of generators of $G$ then 
$|\Sigma(a_1,c_1) \cap \Sigma(a_2,c_2)|\ge
2$.
\end{prop}
Again the 
proof is analogous to the first Proposition of this section.

\begin{footnotesize}
\noindent

\end{footnotesize}

\vfill

\noindent {\bf Authors' 
addresses:}

\bigskip

\noindent Pr. D. Ingrid Bauer and  Prof. Dr. 
Fabrizio Catanese \\
Lehrstuhl Mathematik VIII\\
Universit\"at 
Bayreuth, NWII\\
  D-95440 Bayreuth, Germany

e-mail: 
\\
Ingrid.Bauer@uni-bayreuth.de\\
\indent 
Fabrizio.Catanese@uni-bayreuth.de

\noindent   Prof. Dr. Fritz 
Grunewald \\
Mathematisches Institut 
der\\
Heinrich-Heine-Universit\"at D\"usseldorf\\
 
Universit\"atsstr. 1
                          D-40225 D\"usseldorf, 
Germany\\
e-mail: 
\\grunewald@math.uni-duesseldorf.de

\end{document}